\documentclass[10pt,leqno]{article}
 \usepackage{fancybox} 
 \usepackage{amsmath,latexsym} 
 \usepackage{amscd,amsthm,amsmath,amssymb,amsfonts,times,ascmac}
\usepackage{mathrsfs}
\usepackage[dvips]{color} 
\usepackage[all]{xy}
\usepackage{titlefoot}
\numberwithin{equation}{section}
\theoremstyle{plain}
    \newtheorem{thm}{Theorem}[section]
    \newtheorem{lem}[thm]{Lemma}
    
    \newtheorem{prop}[thm]{Proposition}
    
    \newtheorem{prob}[thm]{Problem}
    \newtheorem{hypo}[thm]{Hypothesis}

    \newtheorem{rem}[thm]{Remark}
  \theoremstyle{definition}
    \newtheorem{defn}[thm]{Definition}
    \newtheorem{exmp}[thm]{Example}
\begin{document}

\def\Tei{{\mathrm{Tei}}}
\def\Coker{\mathrm{Coker}}
\def\crys{\mathrm{crys}}
\def\dlog{{\mathrm{dlog}}}
\def\dR{{\mathrm{d\hspace{-0.2pt}R}}}            
\def\et{{\mathrm{\acute{e}t}}}  
\def\Frac{{\mathrm{Frac}}}
\def\id{{\mathrm{id}}}              
\def\Image{{\mathrm{Im}}}        
\def\Hom{{\mathrm{Hom}}}  
\def\Ext{{\mathrm{Ext}}}
\def\MHS{{\mathrm{MHS}}}  
\def\rig{{\mathrm{rig}}}  
\def\VMHS{{\mathrm{VMHS}}}  
  
\def\ker{{\mathrm{Ker}}}          
\def\Pic{{\mathrm{Pic}}}
\def\CH{{\mathrm{CH}}}
\def\NS{{\mathrm{NS}}}
\def\Sym{{\mathrm{Sym}}}
\def\End{{\mathrm{End}}}
\def\pr{{\mathrm{pr}}}
\def\Proj{{\mathrm{Proj}}}
\def\ord{{\mathrm{ord}}}
\def\zar{{\mathrm{zar}}}
\def\reg{{\mathrm{reg}}}          %
\def\res{{\mathrm{res}}}          %
\def\Res{\mathrm{Res}}
\def\Spec{\operatorname{Spec}}     
\def\syn{{\mathrm{syn}}}
\def\cont{{\mathrm{cont}}}
\def\ln{{\mathrm{ln}}}
\def\bA{{\mathbb A}}
\def\C{{\mathbb C}}
\def\G{{\mathbb G}}
\def\bE{{\mathbb E}}
\def\bF{{\mathbb F}}
\def\F{{\mathbb F}}
\def\bG{{\mathbb G}}
\def\bH{{\mathbb H}}
\def\bJ{{\mathbb J}}
\def\bL{{\mathbb L}}
\def\cL{{\mathscr L}}
\def\bN{{\mathbb N}}
\def\P{{\mathbb P}}
\def\Q{{\mathbb Q}}
\def\R{{\mathbb R}}
\def\Z{{\mathbb Z}}
\def\cH{{\mathscr H}}
\def\cD{{\mathscr D}}
\def\cD{{\mathscr D}}
\def\cF{{\mathscr F}}
\def\cO{{\mathscr O}}
\def\O{{\mathscr O}}
\def\cR{{\mathscr R}}
\def\cS{{\mathscr S}}
\def\cX{{\mathscr X}}
\def\cY{{\mathscr Y}}
\def\cE{{\mathscr E}}

%
\def\ve{\varepsilon}
\def\vG{\varGamma}
\def\vg{\varGamma}
%
%
%
%
\def\lra{\longrightarrow}
\def\lla{\longleftarrow}
\def\Lra{\Longrightarrow}
\def\hra{\hookrightarrow}
\def\ot{\otimes}
\def\op{\oplus}
\def\wt#1{\widetilde{#1}}
\def\wh#1{\widehat{#1}}
\def\ul#1{\underline{#1}}
\def\ol#1{\overline{#1}}
\def\us#1#2{\underset{#1}{#2}}
\def\os#1#2{\overset{#1}{#2}}

\def\FilFMIC{{\mathrm{Fil}\text{-}F\text{-}\mathrm{MIC}}}
\def\FS{{(F)}}
\def\D{\partial}
\def\whO{\O^\wedge_K}
     
\title{Frobenius structure on hypergeometric equations, $p$-adic polygamma 
values and $p$-adic $L$-values\thanks{This research is supported by KAKENHI JP 23K03025	, 18H05233.}}
\author{Masanori Asakura\footnote{Hokkaido University, Department of Mathematics, 
Sapporo 060-0810, Japan.
e-mail: \texttt{asakura@math.sci.hokudai.ac.jp}}
,\quad Kei Hagihara\footnote{Keio University, 
Department of Mathematics, 
Yokohama, 223-8522, Japan. e-mail: \texttt{kei.hagihara@gmail.com}}}
\date\empty

\maketitle

\begin{abstract}
Recently, Kedlaya proved a formula which explicitly describes the Frobenius structure on certain $p$-adic hypergeometric equations.
In this paper, we present a generalization of his formula, which is applicable to cases the original does not cover.
A striking feature of our generalized formula is that, in these newly covered cases, the Frobenius matrix is expressed by the $p$-adic polygamma values and, consequently, by $p$-adic $L$-values for Dirichlet characters. 
As an application to $p$-adic geometry, we show that, for a projective smooth family whose Picard-Fuchs equation is a hypergeometric one, 
the Frobenius matrix on the corresponding log-crystalline cohomology is described in terms of some values of the logarithmic function and $p$-adic $L$-functions.

\end{abstract}

\unmarkedfntext{2020 \textit{Mathematics Subject Classification}. Primary 12H25; Secondary 11S40, 11S80.\newline\hspace*{2em}\textit{Keywords and phrases}. $p$-adic differential equations, hypergeometric equations, $p$-adic polygamma functions, $p$-adic $L$-functions.}
%
%
%

\section{Introduction}

The mysterious link between zeta values and cohomological invariants, both archimedean and non-archimedean, has long intrigued number theorists, serving as a rich source of inspiration for their work. However, its enigmatic nature continues to defy their efforts to fully comprehend it. In this article, we approach the profound mystery through a promising intermediary: the $p$-adic hypergeometric function.

To clarify our motivation, we begin by setting the situation precisely. Let $F$ be a number field, and $\Z_F\subset F$ the ring of integers. 
Let $\mathcal Q$ be a finite set of primes of $F$ and put $R={\mathcal Q}^{-1}\Z_F$.
 Let $\cX\to \Spec R[[z]]$ be a projective flat $R$-morphism that is smooth over 
 $\Spec R[[z]][z^{-1}]$ and
 such that the central fiber $Y$
 is a simple relative normal crossing divisor
over $R$ with multiplicities lying in $R^\times\cap \Z$.
Then, we have the log smooth morphism $(Y,L_Y)\to (\Spec R,L_R)$
that is the pull-back of the log smooth morphism $(\cX,Y)\to (\Spec R[[z]],(z))$
 (cf. \S \ref{log-crys-sect}),
 \[
 \xymatrix{
 (Y,L_Y)\ar[r]\ar[d]&(\cX,Y)\ar[d]\\
 (\Spec R,L_R)\ar[r]&(\Spec R[[z]],(z)).
 }
 \]
From now on, we write $Y_k:=Y\times_R\Spec k$ for an $R$-ring $k$.
Fix a prime number $p$ and a finite place $v$ of $R$ above $p$, and let $F_v$
(resp. $R_v$) denote the $v$-adic completion of $F$ (resp. $R$) and $k_v$ the residue field.
Under the above setting, the log crystalline cohomology group
\[
H^*_{\text{\rm log-crys}}(Y_{k_v}/R_v)
:=H^*_{\text{\rm log-crys}}((Y_{k_v},L_{Y_{k_v}})/(R_v,L_{R_v}))
\]
is defined (\cite[\S 6]{Ka-log}).
This is the $p$-adic counterpart of the limiting mixed Hodge structures
by Steenbrink \cite{steenbrink}.
If we fix a $p$-th Frobenius $\sigma$ on $(R_v,L_{R_v})$,
then the crystalline Frobenius $\Phi_{\crys,v}$ that is $\sigma$-linear is defined 
in the natural way.

\begin{prob}\label{intro-prob}
Describe explicitly the matrix representation of $\Phi_{\crys,v}$ with respect to a basis
of $H^*_{\text{\rm log-dR}}(Y_F/F)$ under the Hyodo-Kato isomorphism
\begin{equation}\label{Hyodo-Kato}
H^*_{\text{\rm log-crys}}(Y_{k_v}/R_v)\ot\Q
\cong
 F_v\ot_FH^*_{\text{\rm log-dR}}(Y_F/F)
\end{equation}
with the log de Rham cohomology (\cite[Theorem (6.4)]{Ka-log}, 
see also \eqref{comparison2} in below).
\end{prob}

If $Y/R$ is smooth, then one can eliminate ``log'' in the above setting. In this case, various connections with $p$-adic special functions are already known, for example,
Ogus \cite[3.15]{ogus}
shows that the crystalline Frobenius on CM elliptic curves are described by $p$-adic gamma values.

On the other hand, it seems that not many things are known in degenerate case.
In such a case,
the log-crystalline cohomology usually
has mixed weights, and then Problem \ref{intro-prob} is, more or less, concerned with
the $p$-adic Beilinson conjecture by Perrin-Riou.
In particular,
the off-diagonal entries of $\Phi_{\crys,v}$ are expected to be
related to some $p$-adic $L$-values. 
Actually, when $X/S$ is the Dwork pencil of quintic threefolds and $Y$ is the totally degenerate fiber,
Shapiro \cite{shapiro} shows that
the $p$-adic zeta value $\zeta_p(3)$ appears
in the Frobenius matrix.

The purpose of this paper is to give another approach to the study of the relation between Frobenii and $p$-adic $L$-values, through an extension of Kedlaya's result on $p$-adic hypergeometric differential equations.

\medskip

Now we present a brief review of Kedlaya's results, followed by a detailed statement of our own. Let $\cD=\Q_p[z,(z-z^2)^{-1},\frac{d}{dz}]$ be the ring of differential operators.
We introduce the hypergeometric differential operator
\[
P(\ul a;\ul b):=\prod_{i=1}^n(D+b_i-1)-z\prod_{i=1}^n(D+a_i),\quad D:=z\frac{d}{dz}
\] 
for parameters $\ul a=(a_1,\ldots,a_n),\ul b=(b_1,\ldots,b_n)\in \Z_{(p)}^n$ and a left $\cD$-module $H_{P(\ul a;\ul b)}=\cD/\cD P(\ul a;\ul b)$. It has a hypergeometric series
\[
F(\ul a;\ul b;z):=\sum_{i=0}^\infty \frac{(a_1)_i\cdots(a_n)_i}{(b_1)_i\cdots(b_n)_i}z^i
\]
as a solution, where $(\alpha)_i:=\alpha(\alpha+1)\cdots(\alpha+i-1)$ denotes the Pochhammer symbol,
cf. \S \ref{notation.sect}.

Defining $\O^\dag=\Q_p[z,(z-z^2)^{-1}]^\dag$ to be the weak completion of $\O=\Q_p[z,(z-z^2)^{-1}]$, we put
$H^\dag_{P(\ul a;\ul b)}=\O^\dag\ot_\O H_{P(\ul a;\ul b)}$.
Then the {\it Frobenius intertwiner} on a hypergeometric equation
is defined to be an additive homomorphism
\[
\Phi:H^\dag_{P(\ul a^{(1)};\ul b^{(1)})}\lra H^\dag_{P(\ul a;\ul b)}
\]
satisfying certain conditions (see \S \ref{sect.Frob.int} for the detail), where 
$(-)^{(1)}$ denotes the Dwork prime 
(cf. \eqref{Dwork.prime})
and we write $\ul a^{(1)}=(a_1^{(1)},\ldots,a_n^{(1)})$.
Under some mild assumption on $\ul a$ and $\ul b$,
it follows from Dwork's uniqueness \cite{Dwork-unique} that
 the Frobenius intertwiner
is unique up to a constant as long as it exists.
On the other hand, it is another question to show
the existence of the Frobenius intertwiner.
Miyatani \cite{Miy} shows it by the method of
middle convolutions. A more natural approach will be constructing a motive associated 
to a hypergeometric equation and applying the theory of $F$-isocrystals
(e.g. \cite{BCM}, \cite{RR}, see also \cite{As-Ross}), though 
the authors do not catch up all works in this direction.

In his paper \cite{Ke},
Kedlaya constructs the Frobenius intertwiner according to an idea of Dwork.
The advantage of this construction is that one can
obtain an explicit formula which describes 
 the residue of $\Phi$ at $z=0$, which we call {\it Kedlaya's residue formula}
 (\cite[Theorem 4.3.2, Corollary 4.3.3]{Ke}, see also Theorem \ref{Kedlaya} below). 
His formula tells that the residue is described by a product of special values
of the $p$-adic gamma function
\[
\Gamma_p(z):=\lim_{n\in\Z_{>0},n\to z}(-1)^n\prod_{0<k<n,\,p\nmid k}k,\quad(z\in\Z_p).
\]
Kedlaya also gives an alternative computation 
for Shapiro's result
in \cite[Example 4.4.1]{Ke}.

However, some cases,
such as the hypergeometric equations of
$F(a_1,...,a_n;1,...,1;z)$, are not covered in his formula.
Such cases often appear in geometric situation,
for example, the connection on the Rham cohomology of the Dwork pencil
\[
x_1^{n}+\cdots+x_{n}^{n}-n\lambda x_1\cdots x_n=0,\quad n\geq 3
\]
of Calabi-Yau hypersurfaces
is the hypergeometric equation of $F(\frac{1}{n},\frac{2}{n},\ldots,\frac{n-1}{n};1,\ldots,1;z)$
(e.g. \cite[Appendix I]{Katz}).
In this paper, we give a generalization of Kedlaya's residue formula that overcomes
these defects.
The precise statement of our formula is provided in Theorem \ref{main.formula1}, together with
a slight generalization Theorem \ref{main-c}.
A distinguished feature is that the residue  
is no longer described by the $p$-adic gamma function, but by {\it$p$-adic polygamma} functions which are defined as follows,
\[
\wt{\psi}_p^{(r)}(z):=\lim_{n\in\Z_{>0},n\to z}\sum_{0< k<n,p\nmid  k}
\frac{1}{k^{r+1}},\quad(z\in\Z_p).
\]
We refer to \cite[2.2]{New} for basic properties on $\wt{\psi}_p^{(r)}(z)$.
Since the $p$-adic polygamma values are 
written by linear combinations of special values of the $p$-adic $L$-functions
$
L_p(s,\chi)
$
for Dirichlet characters $\chi$
(see Theorem \ref{zeta-value3}), we can describe the residue formula in terms of them. As a result, we obtain the following connection between Frobenii and $p$-adic $L$-values:

\begin{thm}[cf. Theorem \ref{degenerate.cor}] \label{intro-thm}
Let $X\to S=\Spec F[z,(z-z^2)^{-1}]$ be a projective smooth $F$-morphism
which has an extension $\cX\to\Spec R[[z]]$ at $z=0$ which fits into the above setting.
Let $M\subset H^{n-1}_\dR(X/S)$ be a sub $\cD_F=F[z,(z-z^2)^{-1},\frac{d}{dz}]$-module.
Suppose that there is an isomorphism
\[
M\cong \bigoplus_{\ul a\in I} H_{P(\ul a;\ul 1)}
\]
where $\ul1:=(1,\ldots,1)$
and $I$ is a finite set of parameters $\ul a=(a_1,\ldots,a_n)$ such that
$0<a_i<1$ and $a_i\in \frac1N\Z$ with $p\nmid N$.
Let $M_{\text{\rm log-crys},v}\subset H^{n-1}_{\text{\rm log-crys}}(Y_{k_v}/R_v)\ot\Q$
be the corresponding subspace defined in \eqref{log-crys-eq3} below. 
Let $\sigma$ be the $p$-th Frobenius on $(R_v,L_{R_v})$ that
is induced from $\sigma(z)=z^p$.
Then, under the Hyodo-Kato isomorphism \eqref{Hyodo-Kato},
the matrix representation of $\Phi_{\crys,v}$ on $M_{\text{\rm log-crys},v}$
with respect to a basis of $H^{n-1}_{\text{\rm log-dR}}(Y_F/F)$ is described by
elements of an algebra generated by some algebraic numbers and
\[
\log (m^{p-1}),\quad L_p(r,\chi\omega^{1-r}),
\]
where $m$ and $r$ run over positive integers such that
$m|N$ and $2\leq r\leq n-1$, 
and $\chi$ runs over primitive characters of conductor $f_\chi$ such that
$f_\chi|N$.
\end{thm}
In fact, one can derive more explicit description of $\Phi_{\crys,v}$ from
Theorem \ref{zeta-value3} and Theorem \ref{degenerate.thm} \eqref{degenerate.thm.eq3.1}.

\smallskip

Recently, Beukers and Vlasenko \cite{BV} extend
Shapiro's result
to arbitrary dimensional Dwork pencils of Calabi-Yau as
a consequence of the study of
the Frobenius structure on the hypergeometric equation
for $F(\frac{1}{n},\ldots,\frac{n-1}{n};1,\ldots,1;z)$. Our formula provides us with a further generalization of it.

The Hodge theoretic counterpart of Theorem \ref{intro-thm} is obtained by
Kerr \cite{Ker}.
He studies limiting mixed Hodge structures for variations of Hodge structures 
with conifold monodromy etc.,
and relates the periods of them to the motivic gamma functions
by Bloch-Vlasenko.
His result can be applied to 
hypergeometric equations $P(\ul a;\ul1)$ such that $\ul a$ are Galois conjugate, and 
provides the counterpart of Theorem \ref{intro-thm}.

\medskip

This paper is organized as follows.
In \S \ref{polygamma-sect}, we review
$p$-adic polygamma functions from \cite[\S 2]{New}.
In particular, we show that $p$-adic polygamma values are linear combinations
of $p$-adic $L$-values because
it is not written in \cite{New}.
In \S \ref{Kedlaya.sect}, we reformulate
Kedlaya's residue formula. This is a preparation for later sections.
In \S \ref{main.theorem.sect}, we formulate and prove our generalization of Kedlaya's
residue formula (Theorems \ref{main.formula1}
and \ref{main-c}).
In \S \ref{application.sect}, we discuss the Frobenius structure of
the log-crystalline cohomology of a family 
such that its Picard-Fuchs equation is a hypergeometric equation.
In \S \ref{syn.sect}, 
we give another application of our formula to the computation of
syntomic regulators of some explicit Milnor symbols.

\bigskip

\noindent{\textbf{Acknowledgement}}.
The authors thank Matt Kerr for letting us know his paper \cite{Ker}
concerning Theorem \ref{intro-thm}.
They also thank the referee for
 reading this paper carefully, and giving us interesting comments.
 
\medskip

\noindent{\textbf{Convention}}.
As usual, we define
$\sum_{i\in \emptyset}a_i=0$ and $\prod_{i\in\emptyset}a_i=1$
for the empty index set.

\section{$p$-adic gamma, $p$-adic polygamma 
and $p$-adic $L$-functions}\label{polygamma-sect}
Throughout this section, we fix a prime number $p\geq2$.
The {\it $p$-adic gamma function} $\Gamma_p(z)$ is the continuous
function on $\Z_p$ defined by
\[
\Gamma_p(z):=\lim_{n\in\Z_{>0},n\to z}(-1)^n\prod_{0<k<n,\,p\nmid k}k.
\]
See e.g. \cite[Chapter II, 6]{Ko}, \cite[11.6]{Cohen2} for the detail.
For an integer $r$, we call a $p$-adically continuous function
\begin{equation}\label{wt-polygamma-def}
\wt{\psi}_p^{(r)}(z):=\lim_{n\in\Z_{>0},n\to z}\sum_{0< k<n,p\nmid  k}
\frac{1}{k^{r+1}},\quad(z\in\Z_p)
\end{equation}
the {\it $p$-adic polygamma function} (\cite[2.2]{New}).
This is an additive analogue of the $p$-adic gamma function.
We notice that this is different from Diamond's polygamma functions \cite{Diamond}.
By definition, one easily shows the following (\cite[Theorem 2.6]{New}),
\[
\wt\psi^{(r)}_p(0)=\wt\psi^{(r)}_p(1)=0,\quad
\wt\psi^{(r)}_p(z)=(-1)^r\wt\psi^{(r)}_p(1-z),\quad
\wt\psi^{(r)}_p(z+1)-\wt\psi^{(r)}_p(z)=\begin{cases}
z^{-r-1}&z\in \Z_p^\times\\
0&z\in p\Z_p.
\end{cases}
\]
\begin{rem}
In \cite[\S 2]{New}, the functions $\psi^{(r)}_p(z)$ are also introduced,
that differ from $\wt\psi^{(r)}_p(a)$
by constant.
However they are not used in this paper.
\end{rem}

The following describes the relation with the $p$-adic gamma.
\begin{thm}[{\cite[Theorem 2.10]{New}}]
\label{log.beta}
Let $B_p(x,y)=\Gamma_p(x)\Gamma_p(y)/\Gamma_p(x+y)$.
Then
\[
\log B_p(z,\ve)=\sum_{i=1}^\infty\wt\psi^{(i-1)}_p(z)\frac{(-1)^{i}\ve^i}{i}
\]
for $z\in \Z_p$ and $\ve\in p\Z_p$,
where
$\log:1+p\Z_p\to p\Z_p$ is the Iwasawa logarithm.
\end{thm}
\begin{lem}\label{polygamma}
There exists a power series expansion
\[
\wt\psi_p^{(r)}(z)=\sum_{k=1}^\infty \alpha_kz^k
\] 
that converges on $|z|_p<\ve$ with some small $\ve>0$.
\end{lem}
\begin{proof}
This follows from Theorem \ref{log.beta} and the fact that
$\Gamma_p(z)$ has a power series expansion which converges
on $|z|_p\leq p^{-(2p-1)/(p(p-1))}$
(\cite[Theorem 11.6.18, (1)]{Cohen2}).
\end{proof}

Let $L_p(s,\chi)$ denote the $p$-adic $L$-function for a Dirichlet character $\chi$.
We then show that the $p$-adic polygamma values are certain linear combinations
of $p$-adic $L$-values, and its converse formula.

\begin{thm}\label{zeta-value}
Let $\omega$ be the Teichm\"uller character.
Let $N>1$ be an integer prime to $p$.
Then,  
\begin{align*}
&\sum_{k=1}^{N-1}\wt\psi^{(r-1)}_p(k/N)
=
(N-N^r)L_p(r,\omega^{1-r})\,
\big(=(N-N^r)(1-p^{-r})\zeta_p(r)\big),& r\ne1,\\
&\sum_{k=1}^{N-1}\wt\psi^{(0)}_p(k/N)
=-N\, p^{-1}\log(N^{p-1}).
\end{align*}
\end{thm}
\begin{proof}
By \cite[Theorem 2.7]{New}, we have
\[
\sum_{k=1}^{N-1}\wt\psi^{(r-1)}_p(k/N)
=N^{r-1}\sum_{k=1}^{N-1}\sum_{\zeta\in\mu_N\setminus\{1\}}(1-\zeta^{-k})\ln_r^{(p)}(\zeta)
=N^{r}\sum_{\zeta\in\mu_N\setminus\{1\}}\ln_r^{(p)}(\zeta)
\]
for arbitrary $r\in\Z$, where $\mu_N$ is the group of $N$-th roots of unity in $\overline\Q_p^\times$.
Then the required equalities follow from \cite[Lemma 2.1, Theorem 2.4]{New}.
\end{proof}
\begin{thm}\label{zeta-value2}
Let $\chi:\Z\to \overline\Q_p$ be a primitive character 
of conductor $f>1$ prime to $p$.
Then, for any $N>0$ prime to $p$ and divisible by $f$,
\[
L_p(r,\chi\omega^{1-r})=-\frac{1}{N^r}\sum_{k=1}^{N-1}
\chi(k)\,\wt\psi^{(r-1)}_p(k/N),\quad(r\in\Z).
\]
\end{thm}
\begin{proof}
For $\zeta\in\mu_{f}$, let
$g(\chi,\zeta):=\sum_{k=1}^{f}\chi(k)\zeta^{-k}$ be the Gauss sum.
By \cite[Theorem 2.7]{New}, we have
\[
\sum_{k=1}^{f-1}\chi(k)\wt\psi^{(r-1)}_p(k/f)
=f^{r-1}\sum_{k=1}^{f-1}\chi(k)\sum_{\zeta\in\mu_{f}
\setminus\{1\}}(1-\zeta^{-k})\ln_r^{(p)}(\zeta)
=-f^{r-1}\sum_{\zeta\in\mu_{f}\setminus\{1\}}g(\chi,\zeta)\ln_r^{(p)}(\zeta).
\]
Now, apply Coleman's formula \cite[I (3)]{C-dlog}
\[
L_p(r,\chi\omega^{1-r})=\frac{1}{f}\displaystyle\sum_{\zeta\in\mu_{f}\setminus\{1\}}g(\chi,\zeta)\ln_{r}^{(p)}(\zeta),
\]
and then we have the required equality in case $N=f$.
Let $N=fm$ with $p\nmid m$.
Recall
\[
\frac{1}{m^{r}}\sum_{i=0}^{m-1}\wt\psi^{(r-1)}_p(z+\frac{i}{m})
=\wt\psi^{(r-1)}_p(mz)
+\sum_{\ve\in\mu_N\setminus\{1\}}\ln^{(p)}_{r}(\ve)
\]
from \cite[(2.17)]{New}. We then have
\begin{align*}
-\frac{1}{N^r}\sum_{k=1}^{N-1}
\chi(k)\,\wt\psi^{(r-1)}_p(k/N)
&=
-\frac{1}{N^r}\sum_{a=1}^{f-1}\sum_{b=0}^{m-1}
\chi(a+bf)\,\wt\psi^{(r-1)}_p\left(\frac{a+bf}{fm}\right)\\
&=
-\frac{1}{N^r}\sum_{a=1}^{f-1}\chi(a)\sum_{b=0}^{m-1}
\wt\psi^{(r-1)}_p\left(\frac{a}{fm}+\frac{b}{m}\right)\\
&=
-\frac{1}{N^r}\sum_{a=1}^{f-1}\chi(a)\cdot m^r
\wt\psi^{(r-1)}_p\left(\frac{a}{f}\right)\\
&=L_p(r,\chi\omega^{1-r})
\end{align*}
as required.
\end{proof}

\begin{thm}\label{zeta-value3}
Let $N>1$ be an integer prime to $p$, and
$N=l_1^{k_1}\cdots l_n^{k_n}$ the prime decomposition.
For a primitive character $\chi$ of conductor $f_\chi$ such that 
$f_\chi\mid N$, put
\[
L_r(\chi)=
\begin{cases}
\prod_i\left(1-\chi(l_i)l_i^{-r}\right)\, L_p(r,\chi\omega^{1-r})& f_\chi>1\\
\left(\prod_i(1-l^{-r}_i)-N^{1-r}\prod_i(1-l_i^{-1})
\right)L_p(r,\omega^{1-r})& f_\chi=1,\, r\ne1\\
\prod_i(1-l_i^{-1})
\left(
\sum_i\frac{p^{-1}\log(l^{p-1}_i)}{l_i-1}+p^{-1}\log(N^{p-1})
\right)
& f_\chi=1,\, r=1.
\end{cases}
\]
Then, for $0<k<N$ with $\gcd(k,N)=1$,
we have
\begin{equation}\label{zeta-value3-eq}
\wt\psi^{(r-1)}_p(k/N)=-\frac{N^r}{\varphi(N)}\sum_{\chi}\chi(k)^{-1}L_r(\chi)
\end{equation}
where $\chi$ runs over all primitive characters of conductor $f_\chi$ such that 
$f_\chi\mid N$.
\end{thm}
\begin{proof}
For a subset $I\subset\{l_1,l_2,\ldots,l_n\}$, put $l_I:=\prod_{l\in I}l$ and 
$N_I:=N/l_I$ (note $l_\emptyset=1$ by convention).
Then, since
\[
\sum_{k=1}^NA_k=\sum_{\gcd(k,N)=1}A_k+\sum_{i}\sum_{l_i|k}A_k
-\sum_{i<j}\sum_{l_il_j|k}A_k
+\cdots
\iff
\sum_{\gcd(k,N)=1}A_k=
\sum_{i=0}^n(-1)^i\sum_{\sharp I=i}\sum_{l_I|k}A_k,
\]
one has
\begin{align*}
\sum_{\gcd(k,N)=1}
\chi(k)\,\wt\psi^{(r-1)}_p(k/N)
=\sum_{i=0}^n(-1)^{i}\sum_{\sharp I=i}
\sum_{l_I|k}
\chi(k)\,\wt\psi^{(r-1)}_p(k/N)
=
\sum_{i=0}^n(-1)^{i}\sum_{\sharp I=i}
\chi(l_I)\sum_{k=1}^{N_I}
\chi(k)\,\wt\psi^{(r-1)}_p(k/N_I).
\end{align*}
Suppose $f_\chi>1$.
If $f_\chi\mid N_I$, then
\[
\sum_{k=1}^{N_I}
\chi(k)\,\wt\psi^{(r-1)}_p(k/N_I)=-N_I^rL_p(r,\chi\omega^{1-r})
\]
by Theorem \ref{zeta-value2} (note $\wt\psi_p^{(r-1)}(1)=0$), and hence
\[
\chi(l_I)\sum_{k=1}^{N_I}
\chi(k)\,\wt\psi^{(r-1)}_p(k/N_I)=
-\chi(l_I)N_I^rL_p(r,\chi\omega^{1-r})=
-\frac{\chi(l_I)}{l_I^r}N^rL_p(r,\chi\omega^{1-r}).
\]
This is also true when $f_\chi\nmid N_I$ as
$\chi(l_I)=0$. Therefore we have
\begin{equation}\label{zeta-value3-eq1}
\sum_{\gcd(k,N)=1}
\chi(k)\,\wt\psi^{(r-1)}_p(k/N)
=
-\sum_{i=0}^n\sum_{\sharp I=i}(-1)^i\frac{\chi(l_I)}{l_I^r}N^rL_p(r,\chi\omega^{1-r})
=
N^r\prod_{i=1}^n(1-\frac{\chi(l_i)}{l_i})\cdot L_p(r,\chi\omega^{1-r}).
\end{equation}
Suppose $f_\chi=1$.
Put
\[
Q_{N_I}=\begin{cases}(N_I-N_I^r)L_p(r,\omega^{1-r})&r\ne 1\\
-N_I\, p^{-1}\log(N_I^{p-1})&r=1.
\end{cases}
\]
In the same way as before,
one has
\begin{equation}\label{zeta-value3-eq2}
\sum_{\gcd(k,N)=1}
\wt\psi^{(r-1)}_p(k/N)
=
\sum_{i=0}^n \sum_{\sharp I=i}(-1)^i\sum_{k=1}^{N_I}
\wt\psi^{(r-1)}_p(k/N_I)
=
\sum_{i=0}^n \sum_{\sharp I=i}(-1)^iQ_{N_I}
\end{equation}
by Theorem \ref{zeta-value}.
If $r\ne1$, then the last term in \eqref{zeta-value3-eq2} is
\[\left(N\prod_i(1-\frac{1}{l_i})-N^r\prod_i(1-\frac{1}{l^r_i})
\right)L_p(r,\omega^{1-r}),\]
and if $r=1$, then
\[
N^r\prod_i(1-\frac{1}{l_i})
\left(
\sum_i\frac{p^{-1}\log(l^{p-1}_i)}{l_i-1}+p^{-1}\log(N^{p-1})
\right).
\]
Now \eqref{zeta-value3-eq} can be derived
from \eqref{zeta-value3-eq1} and \eqref{zeta-value3-eq2}
by the orthogonal relations for characters.
\end{proof}

\section{Kedlaya's residue formula revisited}\label{Kedlaya.sect}
In this section, 
we reformulate the residue formula 
provided in \cite[Theorem 4.3.2, Corollary 4.3.3]{Ke}.
Unfortunately, the formula is not written explicitly in his paper.
We introduce
a basis $\{\wh\omega_k(\ul a;\ul b)\}$ in \S \ref{wh.omega.sect}, and 
reformulate the residue formula explicitly in terms of the basis.
See Remark \ref{hagihara} for the comparison of the original statement by Kedlaya to ours.
This is a preparation for our generalization of the residue formula and its application.
We also give a precise construction of the Frobenius intertwiners
in terms of arithmetic 
differential operators, because details are skipped in \cite{Ke}.

\subsection{Notation}\label{notation.sect}
We fix a prime $p\geq 2$.
Let
$\cD=\cD_{\Q_p}=\Q_p[z,(z-z^2)^{-1},\frac{d}{dz}]$ be 
the ring of differential operators on $\O:=\Q_p[z,(z-z^2)^{-1}]$.
We often use notation
$\partial=\frac{d}{dz}$ and $D=z\partial$.
The action of $\cD$ on a rational function is denoted by $\bullet$,
\[
\cD\times\Q_p(z)\lra \Q_p(z),\quad (P,f)\longmapsto P\bullet f.
\]
For example, 
$Dz^n=z^nD+nz^n$ denotes an element in $\cD$, while
$D\bullet z^n=nz^n$ is a function.
For a field $K\supset \Q_p$,
let $\O_K:=K\ot_{\Q_p}\O$ and $\cD_K:=K\ot_{\Q_p}\cD$ denote the scalar extension.
The action of $\cD_K$ extends 
on $K((z))$, $K((z^{1/n}))$ etc. in a natural way, which we also write
by $\bullet$.

Let $W=W(k)$ be the Witt ring of a perfect field $k$ of characteristic $p>0$,
and $K=\Frac(W)$ the fractional field. 
Let $\mathrm{Frob}_W$ denote the $p$-th Frobenius
on $W$ or $K$.
Let 
\[
\O^\dag_K=K[z,(z-z^2)^{-1}]^\dag=\left\{\sum_{i,j\geq0}c_{ij}z^i(z-z^2)^{-j}
\mid
\text{the radius of convergence of $\sum c_{ij}x^iy^j$ is }>1\right\} 
\] denote the weak completion of $\O_K$ (cf. \cite[p.135]{LS}).
A {\it $p$-th Frobenius} on $\O_K^\dag$ means a $\mathrm{Frob}_W$-linear ring homomorphism
$\sigma:\O_K^\dag\to\O_K^\dag$ 
such that $\sigma(z)-z^p\in pW[z,(z-z^2)^{-1}]^\dag$, where 
``$\mathrm{Frob}_W$-linear'' means $\sigma(\alpha f)=\mathrm{Frob}_W(\alpha)\sigma(f)$ for $\alpha\in K$,
$f\in \O_K^\dag$.
For $P\in \cD_K$, we write
\[
H_P:=\cD_K/\cD_K P,\quad 
H^\dag_P:=\O_K^\dag\ot_{\O_K}H_P.
\]
For $m\in \Z_{\geq0}$, let $(\alpha)_m$ denote the Pochhammer symbol,
that is defined by $(\alpha)_m=\prod_{i=1}^m(\alpha+i-1)$ for $m\geq 1$
and $(\alpha)_0=1$.
For parameters
$\ul{\alpha}=(\alpha_1,\ldots,\alpha_n)$ and
$\ul{\beta}=(\beta_1,\ldots,\beta_n)$ where $\alpha_i,\beta_j\in \C_p$ with $\beta_j\not\in\Z_{\leq 0}$,
let
\[
F(\ul{\alpha};\ul{\beta};z)=\sum_{i=0}^\infty
\frac{(\alpha_1)_i\cdots(\alpha_n)_i}{(\beta_1)_i\cdots(\beta_{n})_i}z^i
\]
denote
the hypergeometric series.
We write
\begin{equation}\label{HG.diff.op}
P(\ul{\alpha};\ul{\beta})=\prod_{j=1}^{n}(D+\beta_j-1)-z\prod_{i=1}^{n}(D+\alpha_i)\in\C_p\ot_{\Q_p}\cD
\end{equation}
for the parameters $\ul \alpha,\ul \beta$,
and call it the {\it hypergeometric differential operator}.
It is immediate to show
\begin{equation}\label{hyp-nonzero.const}
P(\ul{\alpha};\ul{\beta})\bullet F(\ul{\alpha};\ul{\beta};z)=\prod_{j=1}^n(\beta_j-1).
\end{equation}
In particular, if $\beta_j=1$ for some $j$, then $P(\ul{\alpha};\ul{\beta})$ gives the
differential equation of $F(\ul{\alpha};\ul{\beta};z)$ (cf. \cite[2.1]{bailey}).

\subsection{Frobenius intertwiners}\label{sect.Frob.int}
Let $\tau$ be a $p$-th Frobenius on $\O_K^\dag$.
Let $P,Q\in \cD_K$ satisfy that $H_P$ and $H_Q$ are locally free
of finite rank over $\O_K$.
We call an additive homomorphism
\[
\xymatrix{
\Phi:H^\dag_{Q}\ar[r]&
H^\dag_P}
\]
a {\it $\tau$-linear Frobenius intertwiner}  if the following conditions are satisfied.
\begin{itemize}
\item[(i)]
$\Phi$ is $\tau$-linear, which means that $\Phi(fv)=\sigma(f)\Phi(v)$
for $f\in \O_K^\dag$, $v\in H^\dag_Q$,
\item[(ii)]
$D\Phi=p\Phi D$,
\item[(iii)]
$\Phi$ induces a bijection $\O^\dag_K\ot_{\tau}H^\dag_Q\os{\cong}{\to}
H^\dag_P$.
\end{itemize}
Thanks to Dwork's uniqueness theorem \cite{Dwork-unique}, if
$H_P$ is an irreducible $\cD_K$-module with regular singularities
and has exponents in $\Z_{(p)}=\Z_p\cap\Q$,
then the Frobenius 
intertwiner is unique up to multiplication by a nonzero constant as long as it exists.
The hypergeometric differential equation $P(\ul \alpha;\ul \beta)$
 has regular singularities at $z=0,1,\infty$.
If $\ul \alpha,\ul\beta\in \Z_{(p)}^n$, then it has exponents
in $\Z_{(p)}$.
Moreover, thanks to the main theorem of \cite{BH} (see also \cite[Proposition 3.1.11]{Ke}),
$H_{P(\ul \alpha;\ul \beta)}$ is irreducible if $\alpha_i-\beta_j\not\in\Z$.
Therefore, the Frobenius intertwiner on $P(\ul\alpha;\ul\beta)$ is unique
up to $K^\times$ as long as it exists.

\subsection{Construction of the Frobenius intertwiners on hypergeometric equations}\label{Frobenius.sect}
We fix a $p$-th Frobenius $\sigma$ on $\O^\dag=\Q_p[z,(z-z^2)^{-1}]^\dag$
given by $\sigma(z)=z^p$.
We review from \cite[\S 3]{Ke} the construction of
the $\sigma$-linear Frobenius intertwiner on hypergeometric equations.
For $\alpha\in\Z_p$, let $T_\alpha$ denote the unique integer such that
$\alpha+T_\alpha\equiv0$ mod $p$ and $0\leq T_\alpha<p$.
We define
\begin{equation}\label{Dwork.prime}
\alpha^{(1)}=(\alpha+T_\alpha)/p,\quad \alpha^{(i)}=(\alpha^{(i-1)})^{(1)},\quad i\geq 2
\end{equation}
and call $\alpha^{(i)}$ the $i$-th Dwork prime.
For $\ul \alpha,\ul\beta\in \Z_p^{n}$,
we write $\ul\alpha^{(i)}=(\alpha_1^{(i)},\ldots,\alpha_n^{(i)})$ and $\ul\beta^{(i)}
=(\beta_1^{(i)},\ldots,\beta_n^{(i)})$.

\medskip

 Let $\{A_k\}$ be defined by a series expansion
\[
\exp\left(z+\frac{z^p}{p}\right)=\sum_{k=0}^\infty A_k\, z^k.
\]
Let $\pi\in \ol \Q_p$ be an element such that $\pi^{p-1}=-p$.
Substituting $z\mapsto \pi z$, the above becomes 
the Dwork exponential $\exp(\pi(z-z^p))$,
and it is known that it has radius of convergence $p^{(p-1)/p^2}$ (e.g. \cite[I, \S 3]{Ko}).
Therefore one has a lower estimate
\begin{equation}\label{ak-eq1}
\liminf_{k\to\infty}\frac{\ord_p(k!A_k)}{k}\geq\frac{p-1}{p^2}.
\end{equation}
For an differential operator $Q$ and an integer $T\geq 0$, we put
\[
\theta(Q)_T:=\sum_{i=0}^\infty A_i\,Q(Q+1)\cdots(Q+i+T-1).
\]
We only apply $Q=D+\alpha$ or $-D+\alpha$ with $\alpha\in\Z_p$.
Then, the operators are given as follows
\begin{align}\label{theta1}
\theta(D+\alpha)_T
&=\sum_{k=0}^\infty
\left(\sum_{i=0}^\infty
(i+k)!A_{i-T+k}\frac{(\alpha+k)_{i}}{i!}\right)\frac{z^k\D^k}{k!}\\
\theta(-D+\alpha)_T
&=\sum_{k=0}^\infty\left(\sum_{i= 0}^\infty
(-1)^k(i+k)!A_{i-T+k}\frac{(\alpha)_i}{i!}\right)
\frac{z^k\D^k}{k!}\label{theta2}
\end{align}
where $A_j:=0$ for $j<0$, so that they belong to the ring
\begin{align*}
\varGamma(\wh\P^1_{\Z_p},\cD^\dag_{\wh\P^1_{\Z_p},\Q}({}^\dag \infty))
=\left\{\sum_{i,j=0}^\infty c_{ij}z^i\frac{\D^j}{j!}\,\bigg|\,
\text{the radius of convergence of $\sum c_{ij}x^iy^j$ is }>1\right\} 
\end{align*}
by \eqref{ak-eq1}, where $\cD^\dag_{\wh\P^1_{\Z_p},\Q}({}^\dag E)$
denotes 
the sheaf of arithmetic differential operators with overconvergent singularity at a divisor $E$
(\cite{Ber}).
\begin{thm}[{\cite[4.1.3]{Miy}}]\label{Miyatani}
Let $
\cD^\dag:=\varGamma(\wh\P^1_{\Z_p},\cD^\dag_{\wh\P^1_{\Z_p},\Q}({}^\dag \{0,1,\infty\}))$.
Let $\alpha_i,\beta_j\in \Z_{(p)}=\Z_p\cap\Q$ and 
suppose that $\alpha_i-\beta_j\not\in\Z$ for every $i,j$.
Then the natural map
\[\xymatrix{
H^\dag_{P(\ul\alpha;\ul\beta)}=\O_K^\dag\ot \cD/\cD P(\ul a;\ul b) \ar[r]&
\cD^\dag/\cD^\dag P(\ul\alpha;\ul\beta)}\] 
is
bijective.
\end{thm}
\paragraph{Construction of Frobenius intertwiner, \cite[Theorem 4.1.2]{Ke}.}
Let $\alpha_i,\beta_j\in \Z_{(p)}$ satisfy $\alpha_i-\beta_j\not\in\Z$ for every $i,j$.
We construct the $\sigma$-linear Frobenius intertwiner
\begin{equation}\label{FI}
\xymatrix{\Phi:H^\dag_{P(\ul\alpha^{(1)};\ul\beta^{(1)})}\ar[r]&
H^\dag_{P(\ul{\alpha};\ul{\beta})}}
\end{equation}
in the following way\footnote{Kedlaya further discusses the Frobenius intertwiners
for the differential operators 
$P(\{\alpha_1,\ldots,\alpha_n\};\{\beta_1,\ldots,\beta_m\})$
with $n\ne m$.
Unfortunately, the authors do not yet catch up it, 
so we only consider the case $n=m$ throughout this paper.
}.
Let
$\omega\in H_{P(\ul{\alpha};\ul{\beta})}$ and 
 $\omega^\FS\in H_{P(\ul\alpha^{(1)};\ul\beta^{(1)})}$
 denote the elements\footnote{
We use suffix ``$\FS$'' to indicate that 
the object is defined from the parameters $\ul \alpha^{(1)},\ul \beta^{(1)}$.}
 represented by $1_\cD\in\cD$.
Put
\begin{equation}\label{Frobenius-structure-1}
\theta(\ul\alpha;\ul\beta):=\prod_{i=1}^n
\theta(D+\alpha_i)_{T_{\alpha_i}}\cdot\prod_{j=1}^n\theta(-D+1-\beta_j)_{T_{1-\beta_j}}
\in\cD^\dag.
\end{equation}
There are overconvergent functions $E_i(z)\in \Q_p[z,(z-z^2)^{-1}]^\dag$ such that
\begin{equation}\label{Frobenius-structure-2}
\theta(\ul\alpha;\ul\beta)\equiv \sum_{i=0}^{n-1}E_i(z)D^i\mod
\cD^\dag P(\ul{\alpha};\ul{\beta})
\end{equation}
by Theorem \ref{Miyatani}.
Let $\sigma'(z)=(-1)^{n(p-1)}z^p$.
Then \cite[Theorem 4.1.2]{Ke} asserts that
a $\sigma'$-linear map
\[
\xymatrix{
\Phi':H^\dag_{P(\ul\alpha^{(1)};\ul\beta^{(1)})}\ar[r]&
H^\dag_{P(\ul{\alpha};\ul{\beta})}}
\]
that satisfies
\[
\Phi'(\omega^\FS)=\sum_{i=0}^{n-1}E_i(z)D^i\omega,\quad D\Phi=p\Phi D
\]
is well-defined. This is the $\sigma'$-linear Frobenius intertwiner.
If $p\ne2$ or $n$ is even, then $\sigma'=\sigma$ and we set $\Phi=\Phi'$
the required Frobenius intertwiner \eqref{FI}.
If $p=2$ and $n$ is odd, then 
$\sigma'(z)=-z^2\equiv \sigma(z)$ mod $2$, so that the well-known transformation formula of Frobenius allows us 
to obtain the $\sigma$-linear Frobenius intertwiner \eqref{FI}
from $\Phi'$ (cf. \S \ref{changing.sect} below).

\medskip

The above construction is based on the ``exponential module structure''
on the GKZ equations (see \cite[3.3]{Ke} for the detail).
Let $W=\Q_p[x_1,\ldots,x_{2n},\partial_1,\ldots,\partial_{2n}]$ $(\partial_i:=\partial/\partial x_i)$
be the Weyl algebra.
We take a matrix $A=\begin{pmatrix}
I_n&O&1\\ O&-I_{n-1}&1\end{pmatrix}$.
Let $\ul\delta=(\delta_1,\ldots,\delta_{2n-1})\in\Q_p^{2n-1}$.
The GKZ equation for $A$ and $\ul\delta$ is the system of differential equations
\begin{equation}\label{GKZ-hyp}
(\partial_1\cdots\partial_n
-\partial_{n+1}\cdots\partial_{2n})y=0,\quad
\begin{cases}
(x_i\D_i+x_{2n}\D_{2n}+\delta_i)y=0&i=1,\ldots,n\\
(-x_j\D_j+x_{2n}\D_{2n}+\delta_j)y=0&j=n+1,\ldots,2n-1.
\end{cases}
\end{equation}
The pull-back 
by the morphism $z\mapsto(x_1,\ldots,x_{2n})=(1,\ldots,1,z)$
gives rise to the hypergeometric equation
\begin{equation}\label{twisted-HG}
\left(D\prod_{j=n+1}^{2n-1}(D+\delta_j)-(-1)^nz\prod_{i=1}^n(D+\delta_i)\right)y=0
\end{equation}
(cf. \cite[Lemma 3.2.1]{Ke}).
Let $J_{A,\ul\delta}$ be the left ideal of $W$ generated by the differential operators
in the GKZ equation.
Let
\[
R_A=\Q_p[X_1,\ldots,X_n,X_{n+1}^{-1},\ldots,X_{2n-1}^{-1},X_1X_2\cdots X_{2n-1}]
\]
and $R_A[x]=R_A[x_1,\ldots,x_{2n}]$.
Let $X^{(j)}$ denote the $j$-th generator of $R_A$.
Put
\[
g_A=\sum_{j=1}^{2n} \pi x_jX^{(j)},
\]
\[
D_{A,\ul\delta,i}=X_i\frac{\partial}{\partial X_i}+X_i\frac{\partial g_A}{\partial X_i}+\delta_i
\in\mathrm{End}(R_A[x]),
\quad(1\leq i\leq 2n-1).
\]
Then there is an isomorphism
\begin{equation}\label{exp-module}
\xymatrix{
R_A[x]/\sum_i\Image(D_{A,\ul\delta,i})
\ar[r]^-\cong& W/J_{A,\ul\delta}
}\end{equation}
of $W$-modules given by $(\pi X^{(j)})^k\mapsto \D^k_j$ (\cite[Lemma 3.3.3]{Ke}),
where the $W$-module structure
on the left is given by the twisted differential operators $\D_{A,i}=
\D_i+\D_i\bullet g_A$.
Let $\varphi$ be the endomorphism of $R_A[x]$ induced from $\varphi(x_i)=x_i^p$,
$\varphi(X_i)=X_i^p$.
Then a homomorphism
\[
u\longmapsto \exp(-g_A+\varphi(g_A))\varphi(u)
\]
induces the ``Frobenius intertwiner''
on the GKZ equation from $(A,\ul\delta)$ to $(A,p\ul\delta)$
(\cite[3.4.7]{Ke}), 
and thereby on the hypergeometric equation \eqref{twisted-HG} from ${\ul\delta}$ to 
${p\ul\delta}$.
Replace $z$ with $(-1)^nz$, then the Frobenius intertwiner becomes $\sigma'$-linear
where $\sigma'(z)=(-1)^{n(p-1)}z^p$.
Our arithmetic differential operator ``$\theta(\ul\alpha;\ul\beta)$''
comes from ``$\exp(-g_A+\varphi(g_A))$''. 
To be more precise, as a formal series expression, we have 
\[
\exp(-g_A+\varphi(g_A))
=\prod_{j=1}^{2n}
\exp(-\pi x_jX^{(j)}+\pi x_j^p(X^{(j)})^p)
=\prod_{j=1}^{2n}
\left(\sum_{k=0}^\infty A_k\,(-\pi x_jX^{(j)})^k \right)\]
for $p>2$.
This is sent to
\[
\prod_{j=1}^{2n}
\left(\sum_{k=0}^\infty A_k(-1)^kx_j^k\partial_j^k \right)=
\prod_{j=1}^{2n}
\left(\sum_{k=0}^\infty A_k(-D_j)(-D_j+1)\cdots(-D_j+k-1) \right)
\]
by \eqref{exp-module} where $D_j=x_j\partial_j$. 
Since 
$D_i=-D_{2n}-p\delta_i$ $(1\leq i\leq n)$ and
$D_j=D_{2n}+p\delta_j$ $(n+1\leq j<2n)$
by \eqref{GKZ-hyp}, and $D_{2n}$ corresponds to $D=z\partial_z$ under the pull-back
by
$z\mapsto(x_1,\ldots,x_{2n})=(1,\ldots,1,z)$,
the above turns out to be an arithmetic differential operator
\[
\Theta=\prod_{i=1}^n
\theta(D+p\delta_i)_{0}\cdot\prod_{j=n+1}^{2n-1}\theta(-D-p\delta_j)_{0}\cdot\theta(-D)_{0}
\]
in our notation\footnote{
In case $p=2$, we need correction on $\Theta$.
More precisely, we need to replace $\theta(Q)_T$ with
$\theta'(Q)_T$ in $\Theta$
where $\theta'(Q)_T=\sum A'_iQ(Q+1)\cdots(Q+i+T-1)$,
and $A'_k\in\Q$ are defined by $\exp(x-x^2/2)=\sum A'_kz^k$.
However it turns out that the operator $\Theta'=\prod\theta'(D+p\delta_i)_0\cdots$ agrees with $\Theta
=\prod \theta(D+p\delta_i)_0\cdots$ up to sign.
}.
This
induces the Frobenius intertwiner on the equation \eqref{twisted-HG}
from $\ul\delta$ to $p\ul\delta$. To make the Frobenius intertwiner from 
$\ul\delta^{(1)}$ to $\ul\delta$, we take the shift ``$x_j^k\partial_j^k\mapsto
x_j^{k+T}\partial_j^{k+T}$'' according to \cite[Construction 3.4.5]{Ke}.

\begin{lem}\label{cont-lem}
The overconvergent functions $E_i(z)$ in \eqref{Frobenius-structure-2} are 
continuous for
$\ul \alpha$ and $\ul\beta$ with respect to the sup norm $||\sum c_i z^i||:=\sup\{|c_i|\}_i$.
\end{lem}
\begin{proof}
Let $C_k(\alpha)$ be the coefficient of $z^k\D^k/k!$ in \eqref{theta1}.
Then they are uniformly
continuous for $\alpha\in\Z_p$ by \eqref{ak-eq1}, in other words,
$\alpha\mapsto \sum_{k\geq0} C_k(\alpha)z^k$
is continuous with respect to the sup norm.
The same thing is true for \eqref{theta2}.
Now the assertion follows since the arrow in Theorem \ref{Miyatani}
is $p$-adically continuous.
\end{proof}
\subsection{Basis $\wh\omega_k(\ul a;\ul b)$}\label{wh.omega.sect}
\begin{lem}\label{horizontal.lem}
Let $\ul\alpha,\ul\beta\in \C_p^n$ be parameters 
such that $2-\beta_j\not\in \Z_{\leq0}$ for all $j$.
Let $\check F(z):=F(\ul{1-\alpha};\ul{2-\beta};z)$
where $\ul{1-\alpha}=(1-\alpha_1,\ldots,1-\alpha_n)$ etc.
Let $q_i\in(1-z)^{-1}\C_p[z]$ be defined by
\begin{equation}\label{horizontal.lem.eq1}
P=P(\ul{\alpha};\ul{\beta})=(1-z)(D^{n}+q_{n-1}D^{n-1}+\cdots+q_1D+q_0).
\end{equation}
Define $y_i$ by
\begin{equation}\label{horizontal.lem.eq2}
y_{n-1}=(1-z)\check F(z),\quad y_i=q_{i+1}y_{n-1}-D\bullet y_{i+1},\quad (0\leq i\leq n-2)
\end{equation}
Let $H_P:=(\C_p\ot_{\Q_p}\cD)/(\C_p\ot_{\Q_p}\cD) P$, and 
let
$\omega\in H_P$ be the element represented by $1_\cD$.
Suppose $\beta_j=1$ for some $j$.
Then $y_{n-1}D^{n-1}\omega+\cdots+y_1D\omega+y_0\omega$ is a horizontal section
of $\C_p((z))\ot_{\C_p[z]} H_P$,
\[
D(y_{n-1}D^{n-1}\omega+\cdots+y_1D\omega+y_0\omega)=0.
\]
\end{lem}
\begin{proof}
Since
\begin{align*}
&D(y_{n-1}D^{n-1}\omega+\cdots+y_1D\omega+y_0\omega)\\
=&
y_{n-1}D^{n}\omega+(D\bullet y_{n-1}+y_{n-2})D^{n-1}\omega+\cdots+
(D\bullet y_1+y_0)D\omega+(D\bullet y_0)\omega\\
=&(D\bullet y_{n-1}+y_{n-2}-q_{n-1}y_{n-1})D^{n-1}\omega+\cdots+(D\bullet y_1+y_0-q_1y_{n-1})
D\omega+(D\bullet y_0-q_0y_{n-1})\omega
\end{align*}
in $\C_p((z))\ot_{\C_p[z]} H_P$, 
this vanishes if and only if
\[
y_i=q_{i+1}y_{n-1}-D\bullet y_{i+1},\quad (0\leq i\leq n-2),\quad D\bullet y_0-q_0y_{n-1}=0
\]
\[
\iff y_i=q_{i+1}y_{n-1}-D\bullet y_{i+1},\quad (0\leq i\leq n-2),\quad Q\bullet y_{n-1}=0
\]
where 
\[Q:=D^{n}-D^{n-1} q_{n-1}+\cdots+(-1)^{n-1}D  q_1+(-1)^{n}q_0\in\cD.\]
Let $S_i(\alpha)$ denote the symmetric function defined by
$(x+\alpha_1)\cdots(x+\alpha_{n})=x^{n}+S_{n-1}(\alpha)x^{n-1}+\cdots+S_0(\alpha)$.
Then
$(1-z)q_i=S_i(\beta-1)-zS_i(\alpha)$ for $0\leq i\leq n-1$.
We have
\begin{align*}
Q(1-z)&=
D^{n}(1-z)-\sum_{i=0}^{n-1}(-1)^{n-1-i}\left(S_i(\beta-1)D^i-zS_i(\alpha)(D+1)^i\right)\\
&=
D^{n}-\sum_{i=0}^{n-1}(-1)^{n-1-i}S_i(\beta-1)D^i-z
\left((D+1)^{n}+\sum_{i=0}^{n-1}(-1)^{n-1-i}S_i(\alpha)(D+1)^i\right)\\
&=\prod_{i=1}^{n}(D+1-\beta_i)-z\prod_{i=1}^{n}(D+1-\alpha_i),
\end{align*}
and this is the hypergeometric differential operator of $F(\ul{1-\alpha},\ul{2-\beta};z)$.
Therefore, if we put $y_{n-1}=(1-z)F(\ul{1-\alpha},\ul{2-\beta};z)$,
then it is annihilated by $Q$
(here we use the assumption $\beta_j=1$, see \eqref{hyp-nonzero.const}).
\end{proof}
\paragraph{Definition of $\wh\omega_k(\ul a;\ul b)$.}
Let $
\ul a=(a_1,a_2,\ldots,a_n)\in\Q^n$ and $\ul b=(b_1,b_2,\ldots,b_n)\in\Q^n$
be parameters such that
$b_i-b_j\not\in \Z\setminus\{0\}$ for every $i,j$.
For $1\leq k\leq n$, 
we put
\[
G_k(z)=F(\ul{1-b_k+a},\ul{1-b_k+b};z),
\quad \check G_k(z)=F(\ul{b_k-a},\ul{1+b_k-b};z)
\]
where $\ul{1-b_k+a}=(1-b_k+a_1,1-b_k+a_2,\ldots,1-b_k+a_n)$ etc.,
and
\[
P=P(\ul a;\ul b),\quad Q_k=P(\ul{1-b_k+a};\ul{1-b_k+b}).
\]
Write $H_P=\cD/\cD P$ and $H_{Q_k}=\cD/\cD Q_k$ as before.
Since
$z^{b_k-1}P=Q_kz^{b_k-1}$,
we have an isomorphism
\[
\delta_k: \cD/\cD Q_k\os{\cong}{\lra} \cD/\cD P,\quad 
\theta\longmapsto z^{1-b_k}\theta z^{b_k-1}
\]
of $\Q_p[z,(z-z^2)^{-1}]$-modules.
Let $q^{[k]}_i\in (1-z)^{-1}\Q[z]$ and $y_i^{[k]}\in\Q_p[[z]]$ be defined by
\[
Q_k
=(1-z)(D^{n}+q^{[k]}_{n-1}D^{n-1}+\cdots+q^{[k]}_1D+q^{[k]}_0),
\] 
\[
y_{n-1}^{[k]}
=(1-z)\check G_k(z),\quad y_i^{[k]}=q^{[k]}_{i+1}y^{[k]}_{n-1}-
D\bullet y^{[k]}_{i+1},\quad (0\leq i\leq n-2)
\]
in a similar way to \eqref{horizontal.lem.eq1} and \eqref{horizontal.lem.eq2}.
Then
$y_{n-1}^{[k]}D^{n-1}\omega+\cdots+y^{[k]}_1D\omega+y^{[k]}_0\omega$
is a horizontal section of $\wh H_{Q_k}:=\Q_p((z))\ot_{\Q_p[z]} H_{Q_k}$
by  Lemma \ref{horizontal.lem}.
We define
\begin{align}
\wh\omega_k(\ul a,\ul b)&:=
\delta_k(y_{n-1}^{[k]}D^{n-1}\omega+\cdots+y^{[k]}_1D\omega+y^{[k]}_0\omega)\notag\\
&=
y_{n-1}^{[k]}(D+b_k-1)^{n-1}
\omega+\cdots+y^{[k]}_1(D+b_k-1)\omega+y_0^{[k]}\omega\in \wh H_P,\label{wh.omega.defn}
\end{align}
so that this is annihilated by $\delta_k(D)=D+b_k-1$,
\[
(D+b_k-1)\wh\omega_k(\ul a,\ul b)=0.
\]
\begin{lem}\label{wh.omega.lem}
For $1\leq k\leq n$, we have 
\begin{equation}\label{wh.omega.lem.eq1}
D\wh\omega_k(\ul a,\ul b)=(1-b_k)\wh\omega_k(\ul a,\ul b)\in \wh H_P,
\end{equation}
\begin{equation}\label{wh.omega.lem.eq2}
\wh\omega_k(\ul a,\ul b)\bullet z^{1-b_k} G_k(z)=
\left(\prod_{j\ne k}(b_j-b_k)\right)z^{1-b_k}.
\end{equation}
Suppose that $b_j\ne b_k$ for any $j\ne k$ (hence $b_j-b_k\not\in\Z$).
Then
$\{\wh\omega_1(\ul a,\ul b),\ldots,\wh\omega_n(\ul a,\ul b)\}$ forms a $\Q_p((z))$-basis of $\wh H_P$, and
\begin{equation}\label{wh.omega.lem.eq3}
\wh\omega_k(\ul a,\ul b)\bullet z^{1-b_j} G_j(z)=0\quad(j\ne k).
\end{equation}
Furthermore, if $\eta\in \wh H_P$ satisfies
$D\eta=(1-b_k)\eta$, then $\eta=\wh\omega_k(\ul a,\ul b)$ up to 
multiplication by a constant.
\end{lem}
\begin{proof}
We have shown \eqref{wh.omega.lem.eq1} in the above.
We show \eqref{wh.omega.lem.eq2}.
Write $\wh\omega_k=\wh\omega_k(\ul a,\ul b)$ simply.
We have
\begin{align*}
\wh\omega_k\bullet z^{1-b_k} G_k(z)
&=z^{1-b_k}(y_{n-1}^{[k]}D^{n-1}\omega+\cdots+y^{[k]}_1D\omega+y^{[k]}_0\omega)\bullet G_k(z)\\
&=z^{1-b_k}
\overbrace{(y_{n-1}^{[k]}D^{n-1}\bullet G_k(z)+\cdots+y^{[k]}_1D\bullet G_k(z)+y^{[k]}_0G_k(z))}^{g(z)}.
\end{align*}
Applying $D-1+b_k$, one has $D\bullet g(z)=0$ by \eqref{wh.omega.lem.eq1}.
Therefore, $g(z)=c$ is a constant.
We have
\[
c=g(0)
=y^{[k]}_0G_k(z)|_{z=0}=y^{[k]}_0(0).
\]
Since $y_i^{[k]}(0)=q^{[k]}_{i+1}(0)y^{[k]}_{n-1}(0)=q^{[k]}_{i+1}(0)$ by definition of $y_i^{[k]}$,
we have
\[
c=y_0^{[k]}(0)=q^{[k]}_1(0)=\prod_{j\ne k}(b_j-b_k).
\]
This completes the proof of \eqref{wh.omega.lem.eq2}.

Suppose that $b_j-b_k\not\in\Z$ for any $j\ne k$.
Put
\[
\wh H_{P,\Q_p[[z]]}=\sum_{i=0}^n\Q_p[[z]]D^i\omega\subset \wh H_P
\]
a free $\Q_p[[z]]$-module of rank $n$ that satisfies $\Q_p((z))\ot_{\Q_p[[z]]}
\wh H_{P,\Q_p[[z]]}=\wh H_P$.
It is immediate to see that 
$\wh\omega_k|_{z=0}\in
\wh H_{P,\Q_p[[z]]}/z\wh H_{P,\Q_p[[z]]}$ does not vanish, and satisfies
$(D+b_k-1)(\wh\omega_k|_{z=0})=0$ by \eqref{wh.omega.lem.eq1}.
This implies that
$\wh\omega_1|_{z=0}, \ldots,\wh\omega_n|_{z=0}$ are linearly independent in 
$\wh H_{P,\Q_p[[z]]}/z\wh H_{P,\Q_p[[z]]}$ as $b_j\ne b_k$ for any $j\ne k$.
Therefore
$\{\wh\omega_1, \ldots,\wh\omega_n\}$ forms a $\Q_p[[z]]$-basis
of $\wh H_{P,\Q_p[[z]]}$ by Nakayama's lemma.
We show \eqref{wh.omega.lem.eq3}.
Let $f(z)=\wh\omega_k \bullet z^{1-b_j}G_j(z)$. 
One immediately has $f(z)\in z^{1-b_j}\Q_p[[z]]$.
On the other hand, 
\eqref{wh.omega.lem.eq1} implies
$D\bullet f(z)=(1-b_k)f(z)$, and hence
\[
f(z)=\text{(const.)}z^{1-b_k}.
\]
Since $b_j\not\equiv b_k$ mod $\Z$, this imposes $f(z)=0$.

Let $\eta\in \wh H_P$ satisfies $D\eta=(1-b_k)\eta$.
Let $\eta=\sum_if_i\wh\omega_i$ with $f_i\in \Q_p((z))$.
For $j=1,2,\ldots,n$, one has
\begin{align*}
0&=(D+b_k-1)\eta\bullet z^{1-b_j}G_j(z)
=\sum_{i=1}^n(D+b_k-1)f_i\wh\omega_i\bullet z^{1-b_j}G_j(z)\\
&=(D+b_k-1)\bullet\left(\prod_{l\ne j}(b_l-b_j)z^{1-b_j}f_j\right)
\end{align*}
by \eqref{wh.omega.lem.eq2} and \eqref{wh.omega.lem.eq3}.
Therefore, $z^{1-b_j}f_j$ is annihilated by $D+b_k-1$, and this is equivalent to
that $(D+b_k-b_j)\bullet f_j=0$.
If $j\ne k$, this imposes $f_j=0$ as $b_k-b_j\not\in \Z$.
If $j=k$, it turns out that $f_k$ is a constant.
We thus have $\eta=\text{(const.)}\,\wh\omega_k$ as required.
\end{proof}
\subsection{Kedlaya's residue formula}\label{Frob.wh.omega_k.sect}
Let $\ul a=(a_1,\ldots,a_n)\in\Z_{(p)}^n$ and
$\ul b=(b_1,\ldots,b_n)\in\Z_{(p)}^n$ be parameters satisfying the following hypothesis.

\begin{hypo}\label{Hypothesis}
$0<a_i,b_j<1$ and
$a_i-b_j\not\in\Z$ for every $i,j$, and $b_1,\ldots,b_n$ are pairwise distinct.
\end{hypo}

We keep the notation in the previous section.
We note that $\ul a^{(1)}$ and $\ul b^{(1)}$ also satisfy \textbf{Hypothesis} \ref{Hypothesis}.
Let 
\[
\wh\omega_k:=\wh\omega_k(\ul a,\ul b)\in \wh H_{P(\ul a;\ul b)},\quad
\wh\omega^\FS_k:=\wh\omega_k(\ul a^{(1)},\ul b^{(1)})\in \wh H_{P(\ul a^{(1)};\ul b^{(1)})}
\]
be the elements defined in \eqref{wh.omega.defn}. They form bases by Lemma \ref{wh.omega.lem}.
Put
\[
\mu_k:=T_{1-b_k}=p(1-b^{(1)}_k)-(1-b_k)
,\quad c_k:=\prod_{1\leq i\leq n,\,i\ne k}(b_k-b_i) 
,\quad c^\FS_k:=\prod_{1\leq i\leq n,\,i\ne k}(b^{(1)}_k-b^{(1)}_i)
\]
(see the beginning of \S \ref{Frobenius.sect} for the notation $T_\alpha$).
Let $\sigma$ be the $p$-th Frobenius on $\Q_p((z))$ given by $\sigma(z)=z^p$.
Let 
\[
\xymatrix{
\Phi:\wh H_{P(\ul a^{(1)};\ul b^{(1)})}\ar[r]& \wh H_{P(\ul a;\ul b)}
}\]
be the $\sigma$-linear homomorphism induced from the Frobenius intertwiner in \S \ref{Frobenius.sect}.
Note here that, since each $E_i(z)$ actually belongs to $\mathbb{Q}_p[[z]]$, $\Phi'$ is defined over $\tilde{\mathscr{O}}=\mathscr{O}^\dag\cap\mathbb{Q}_p[[z]]$ in the sense that we have a $\sigma'$-linear map
$\tilde{\Phi}: \tilde{\mathscr{O}}\otimes_{\mathscr{O}} H_{P(\ul{\alpha}^{(1)}; \ul{\beta}^{(1)})} \to \tilde{\mathscr{O}}\otimes_{\mathscr{O}} H_{P(\ul{\alpha}; \ul{\beta})}$
such that $\mathscr{O}^\dag\otimes\tilde{\Phi} = \Phi'$.
Then we put ${\Phi}=\Q_p((z))\ot_{\tilde\O}\tilde\Phi$ in case $p>2$ or $n$ even, and 
otherwise we further apply the transformation formula of Frobenius as in \S \ref{Frobenius.sect}.

We now see\begin{equation}\label{gammak}
\Phi(\wh\omega^\FS_k)=\text{(const.)} \times z^{\mu_k}\wh\omega_k.
\end{equation}
Indeed, one has
\[
D\Phi(\wh\omega^\FS_k)=p\Phi D(\wh\omega^\FS_k)
\os{\eqref{wh.omega.lem.eq1}}{=}
p(1-b^{(1)}_k)\Phi(\wh\omega^\FS_k)=
(\mu_k+1-b_k)\Phi(\wh\omega^\FS_k),
\]
and hence
\[
D(z^{-\mu_k}\Phi(\wh\omega^\FS_k))
=(1-b_k)(z^{-\mu_k}\Phi(\wh\omega^\FS_k)).
\]
This implies that
$z^{-\mu_k}\Phi(\wh\omega^\FS_k)$ 
agrees with $\wh\omega_k$
up to multiplication by a constant
by Lemma \ref{wh.omega.lem}.

\medskip

Kedlaya's residue formula tells that 
the constant in \eqref{gammak} is a product of $p$-adic gamma values.
By Lemma \ref{wh.omega.lem},
it can be reduced to compute the constant appearing in the following equality,
\[
\theta(\ul a;\ul b)\bullet z^{1-b_k}G_k(z)=\text{(const.)}\times z^{p(1-b_k^{(1)})}G_k^\FS(z^p).
\]
This is a straightforward computation
by the fact that
\[
\theta(r_0+r_1 D)_T\bullet z^m=\left(\sum_{i=0}^\infty
A_i\cdot(r_0+r_1 m)_{T+i}\right)z^m,\quad(r_0,\,r_1\in\Z_p)
\]
and a formula (cf. \cite[11.6.15]{Cohen2})
\[
\sum_{i=0}^\infty A_i\cdot(x)_{T+i}=\begin{cases}
-\Gamma_p(1-x)&T+x\equiv0\\
0&T+x\not\equiv0
\end{cases}
\]
for $x\in \Z_p$ and $T\in \{0,1,\ldots,p-1\}$.
In this way, we have the residue formula by Kedlaya.

\begin{thm}[{\cite[4.3.3]{Ke}}]\label{Kedlaya}
Let $\ul a,\ul b\in\Z_{(p)}^n$ be parameters satisfying {\rm \textbf{Hypothesis} \ref{Hypothesis}}.
Let $(x)_+=x$ for $x>0$ and $=1$ otherwise.
Put
\[
K(x):=\prod_{i=1}^n\frac{(a_i-x)_+}{(b_i-x)_+},
\]
\[
Z(x)=\sharp\{i\in\{1,2,\ldots,n\}\mid a_i<x\}-\sharp\{j\in\{1,2,\ldots,n\}\mid b_j<x\}
\]
and let $K^\FS$ and $Z^\FS$ denote the functions defined in the same way
from $\ul a^{(1)},\ul b^{(1)}$.
Put
\begin{equation}\label{Kedlaya-gammak}
\gamma_k:=(-1)^{Z(b_k)}p^{-Z^\FS(b^{(1)}_k)}\frac{K^\FS(b^{(1)}_k)}{K(b_k)}
\prod_{i=1}^n\frac{\Gamma_p(\{b_i-b_k\})}{\Gamma_p(\{a_i-b_k\})}.
\end{equation}
Moreover we put $e=\sum_{i=1}^n T_{a_i}$ if $p>2$ and $e=\sum_{i=1}^na_i+a_i^{(1)}$ if $p=2$.
Then
\begin{equation}\label{Kedlaya-gammak-phi}
\Phi((c_k^\FS)^{-1}\wh\omega^\FS_k)=(-1)^{n+e}p\,\gamma_k \,z^{\mu_k}c_k^{-1}\wh\omega_k.
\end{equation}
(Note that if $p=2$, then $e\in \Z_2$ is not necessarily an integer, but $(-1)^e$ makes sense.)
\end{thm}
\begin{rem}\label{FI-rem}
The Frobenius intertwiner on $P(\ul a;\ul b)$
is determined only up to $\Q_p^\times$.
On the other hand, in Theorem \ref{Kedlaya}, we choose $\Phi$ to be the intertwiner
given by $\theta(\ul a;\ul b)$, so that
the formula \eqref{Kedlaya-gammak-phi} is valid without ambiguity of $\Q_p^\times$.
If we allow to replace $\Phi$ with $c\Phi$ by some $c\in\Q_p^\times$,
one can delete ``$(-1)^{n+e}p$'' in \eqref{Kedlaya-gammak-phi} as it does not depend on $k$.
\end{rem}
\begin{rem}\label{hagihara}
Actually, Kedlaya gives the formula for the bundle
\[
E_{P(\ul{a};\ul{b})} = \bigoplus_{i=1}^n\O e_{i}, \quad (\O:=\mathbb{Q}_p[z, (z-z^2)^{-1}])
\]
which is endowed with the connection
\[
\nabla({\boldsymbol v}) = \frac{dz}{z}\ot(D{\boldsymbol v} + N{\boldsymbol v}),
\]
where ${\boldsymbol v} = \sum_{i=1}^nv_ie_i$ is identified 
with a column vector ${}^t\!(v_1, \ldots, v_n)$ of length $n$, 
and the matrix $N$ is defined to be
\[
\begin{pmatrix} 
  0        & -1     & 0 &  \dots & 0 \\
  0        &  0      &   \ddots     &     & \vdots  \\
  \vdots & \vdots  &  & \ddots  &  0 \\
  0        & 0     & \dots &     0     & -1 \\
  q_{0}    & q_{1} & \dots & q_{n-2}  & q_{n-1}
\end{pmatrix}
\]
with $P(\ul{a};\ul{b}) = 
(1-z)(D^n+ q_{n-1}D^{n-1}+\cdots + q_1D+ q_0)$\,{\rm(cf. \cite[Remark 2.1.2]{Ke})}.
This is a dual connection of our $H_{P(\ul a;\ul b)}$ via
the $\mathscr{O}$-bilinear pairing
\[
\xymatrix{
\langle-, -\rangle: H_{P(\ul{a};\ul{b})} \otimes E_{P(\ul{a};\ul{b})} \ar[r]& \mathscr{O}
}\]
given by $\langle D^{i-1}\omega, e_j\rangle = \delta_j^i$ for $1\le i, j \le n$. 
Let
$
\Phi_E: \sigma^\ast{E}^\dag_{P(\ul{a}^{(1)};\ul{b}^{(1)})} \os{\cong}{\to} {E}^\dag_{P(\ul{a};\ul{b})}$
and 
$\Phi_H: \sigma^\ast H^\dag_{P(\ul{a}^{(1)};\ul{b}^{(1)})}\os{\cong}{\to} H^\dag_{P(\ul{a};\ul{b})}$
be the Frobenius intertwiners respectively.
Let us denote 
$(-)^\vee := \mathrm{Hom}_{\mathscr{O}^\dag}(-, \mathscr{O}^\dag)$.
We then have an isomorphism
\[
\xymatrix{
(\Phi_E^{-1})^\vee:
\sigma^\ast H^\dag_{P(\ul{a}^{(1)};\ul{b}^{(1)})} \cong (\sigma^\ast{E}^\dag_{\ul{a}^{(1)};\ul{b}^{(1)}})^\vee \ar[r]& ({E}^\dag_{\ul{a};\ul{b}})^\vee \cong H^\dag_{P(\ul{a};\ul{b})},}
\]
and this turns out to be a Frobenius intertwiner.
Hence $(\Phi_E^{-1})^\vee$
agrees with $\Phi_H$ up to a scalar multiplication by Dwork's uniqueness theorem.
The reader can now
compare Theorem \ref{Kedlaya} with the formula on $\Phi_E$ in {\rm \cite[4.3.3]{Ke}}.
\end{rem}
\begin{lem}\label{gamma.formula}
Let $\gamma_k$ be the constant \eqref{Kedlaya-gammak}. Then 
\begin{equation}\label{gamma.formula.eq1}
\gamma_k=p^{-1}
\prod_{i=1}^n\frac{(1-b_k+a_i)_{\mu_k}}{\Gamma_p(1-b_k+a_i+\mu_k)}
\prod_{j=1}^n
\frac{\Gamma_p(1-b_k+b_j+\mu_k)}{(1-b_k+b_j)_{\mu_k}},
\end{equation}
where $(\alpha)_i=\alpha(\alpha+1)\cdots(\alpha+i-1)$ denotes
the Pochhammer symbol, and $\mu_k:=p(1-b_k^{(1)})-(1-b_k)\in\{0,1,\ldots,p-1\}$.
\end{lem}
\begin{proof}
We define a function $\chi: \mathbb{R}\rightarrow \{0,1\}$ to be 
$\chi(x)=1$ if $x>0$, and $\chi(x)=0$ otherwise. Then we have
\[
Z(x) = \sum_{i=1}^n \chi(x-a_i) - \sum_{i=1}^n \chi(x-b_i),
\] 
and similarly for $Z^{(F)}$.
Hence, setting
\[A(x) := (-1)^{\chi(b_k-x)}p^{\chi(b_k^{(1)}-x^{(1)})}
\frac{(x-b_k)_+}{(x^{(1)}-b_k^{(1)})_+}\Gamma_p(\{x-b_k\}),
\]
we have
\[
\gamma_k= 
\prod_{i=1}^n\frac{A(b_i)}{A(a_i)}.
\]
On the other hand,  
let $B(x)$ be defined by
\[
\frac{\Gamma_p(x-b_k+1+\mu_k)}{(x-b_k+1)_{\mu_k}}
=\frac{(-1)^{\mu_k}\Gamma_p(x-b_k+1)}{B(x)}
\]
for $x\in\Z_{(p)}$ with $0<x<1$.
We then show
\begin{equation}\label{gamma.formula.eq2}
B(x)=
\begin{cases}
(px^{(1)}-pb_k^{(1)})_+&\text{if }x\not\equiv b_k\mod p\\
1&\text{if }x\equiv b_k\mod p.
\end{cases}
\end{equation}
For $a\in\Z_p$,
let $T_a\in \{0,1,\ldots,p-1\}$ denote the integer such that $a+T_a\equiv 0$ mod $p$.
Put 
\[
l:=T_{x-b_k+1}=\begin{cases}
\mu_k+T_x&\mu_k+T_x<p\\
\mu_k+T_x-p&\mu_k+T_x\geq p.
\end{cases}
\]
Since $\Gamma_p(s+1)=-s\Gamma_p(s)$ for $s\not\in p\Z_p$
and $\Gamma_p(s+1)=-\Gamma_p(s)$ for $s\in p\Z_p$, we have that
$B(x)=x-b_k+1+l$ if $0\leq l\leq \mu_k-1$
and $B(x)=1$ otherwise.
If $x\equiv b_k$ mod $p$, then $l=p-1$ so that $B(x)=1$.
Suppose $x\not\equiv b_k$ mod $p$.
Then $l\leq p-2$.
If $\mu_k+T_x<p$, then 
$l=\mu_k+T_x\geq \mu_k$ so that $B(x)=1$.
In this case, $px^{(1)}-pb_k^{(1)}=x-b_k+1-p+l\leq x-b_k-1<0$.
If $\mu_k+T_x\geq p$, then 
$l=\mu_k+T_x-p\leq \mu_k-1$ so that $B(x)=x-b_k+1+l=px^{(1)}-pb_k^{(1)}>0$.
This completes the proof of \eqref{gamma.formula.eq2}.

If $x\not\equiv b_k$ mod $p$, then
\begin{align*}
\frac{\Gamma_p(x-b_k+1+\mu_k)}{(x-b_k+1)_{\mu_k}}
&=\frac{(-1)^{\mu_k}\Gamma_p(x-b_k+1)}{(px^{(1)}-pb^{(1)}_k)_+}\\
&=(-1)^{\mu_k+\chi(x-b_k)}p^{-\chi(x^{(1)}-b_k^{(1)})}
\frac{(x-b_k)_+\Gamma_p(\{x-b_k\})}{(x^{(1)}-b^{(1)}_k)_+}
\\
&=(-1)^{\mu_k+1-\chi(b_k-x)}p^{-1+\chi(b_k^{(1)}-x^{(1)})}
\frac{(x-b_k)_+\Gamma_p(\{x-b_k\})}{(x^{(1)}-b^{(1)}_k)_+}
\\
&=(-1)^{\mu_k+1}p^{-1}A(x).
\end{align*}
If $x\equiv b_k$ mod $p$, then $p(x^{(1)}-b_k^{(1)})=x-b_k$ and 
\[
\frac{\Gamma_p(x-b_k+1+\mu_k)}{(x-b_k+1)_{\mu_k}}
=(-1)^{\mu_k}\Gamma_p(x-b_k+1)=(-1)^{\mu_k+1+\chi(b_k-x)}\Gamma_p(\{x-b_k\}),
\]
and 
\[A(x) 
= (-1)^{\chi(b_k-x)}p^{\chi(b_k^{(1)}-x^{(1)})}
\frac{(x-b_k)_+}{(x^{(1)}-b_k^{(1)})_+}\Gamma_p(\{x-b_k\})
= (-1)^{\chi(b_k-x)}p^{\chi(b_k-x)+\chi(x-b_k)}
\Gamma_p(\{x-b_k\}).
\]
Summing up the above, we have
\[
\frac{\Gamma_p(x-b_k+1+\mu_k)}{(x-b_k+1)_{\mu_k}}=\begin{cases}
(-1)^{\mu_k+1}p^{-1}A(x)&x\ne b_k\\
(-1)^{\mu_k+1}&x= b_k.
\end{cases}
\]
Therefore, the right hand side of \eqref{gamma.formula.eq1} is
\[
p^{-1}\left(\prod_{i=1}^n(-1)^{\mu_k+1}p^{-1}A(a_i)\right)^{-1}
\left(p
\prod_{j=1}^n(-1)^{\mu_k+1}p^{-1}A(b_j)\right)
=\prod_{i=1}^n\frac{A(b_i)}{A(a_i)}=\gamma_k
\]
as required.
\end{proof}
\begin{lem}\label{gamma.formula.cor}
Let $B_p(x,y)=\Gamma_p(x)\Gamma_p(y)/\Gamma_p(x+y)$.
Suppose $b_k\equiv b_l\equiv 1$ mod $p$.
Put $\ve_k=b_k-1$ and $\ve_l=b_l-1$. 
Then
\[
\frac{\gamma_k}{\gamma_l}
=
\prod_{i=1}^n\frac{\Gamma_p(b_i-\ve_k)}{\Gamma_p(a_i-\ve_k)}
\frac{\Gamma_p(a_i-\ve_l)}{\Gamma_p(b_i-\ve_l)}
=
\prod_{i=1}^n\frac{B_p(b_i,-\ve_l)}{B_p(b_i,-\ve_k)}\frac{B_p(a_i,-\ve_k)}{B_p(a_i,-\ve_l)}.
\]
\end{lem}
\begin{proof}
Since $\mu_k=\mu_l=0$, this is immediate from Lemma \ref{gamma.formula}.
\end{proof}

\section{Generalization of Kedlaya's residue formula}\label{main.theorem.sect}
We keep the notation in \S \ref{notation.sect}.
Moreover we write \[
\wh H_P:=\Q_p((z))\ot_{\Q_p[z]}H_P.
\]
Let $\omega\in H_P$ (or $\wh H_P$) denote the element
represented by $1_\cD$.
For a power series $f_b(z)=\sum c_i(b)z^i$ whose coefficients are $p$-adically 
continuous functions of
variables $b=(b_1,\ldots,b_n)$, we mean by limit $b\to b_0$ the {\it coefficient-wise limit},
\[
\lim_{b\to b_0} f(z)=\sum_i\left(\lim_{b\to b_0} c_i(b)\right)z^i.
\] 
For $\xi=f_0(z)\omega+f_1(z)D\omega+\cdots+f_{n-1}D^{n-1}\omega \in \wh H_P$,
let the limit of $\xi$ be defined by
\[
\lim_{b\to b_0} \xi=(\lim_{b\to b_0} f_0(z))\omega
+(\lim_{b\to b_0} f_1(z))D\omega+\cdots+(\lim_{b\to b_0} f_{n-1}(z))D^{n-1}\omega.
\]
\subsection{Definition of $\wh\omega(m)$}\label{base.sect}
Let $\ul a=(a_1,\ldots,a_n), \ul b=(b_1,\ldots,b_n)\in\Q_p^n$ be parameters such that 
$b_i-b_j\not\in\Z$ for any $i\ne j$.
Let
$P=P(\ul a;\ul b)$ be the hypergeometric differential operator in \eqref{HG.diff.op}.
Let
\[
\wh\omega_k=\wh\omega_k(\ul a,\ul b)=y_{n-1}^{[k]}(D+b_k-1)^{n-1}
\omega+\cdots+y^{[k]}_1(D+b_k-1)\omega+y_0^{[k]}\omega\in \wh H_P
\]
be as in \eqref{wh.omega.defn}.
For a subset $I\subset\{1,2,\ldots,n\}$, we put
\[
\wh\omega_I:=\sum_{k\in I}\left(\prod_{i\in I\setminus\{k\}}\frac{1}{b_k-b_i}\right)\wh\omega_k.
\]
We note that
$\sum_{i\in \emptyset}(-)=0$ and $\prod_{i\in\emptyset}(-)=1$ by convention, 
\[
\wh\omega_\emptyset=0,\quad
\wh\omega_{\{k\}}=\wh\omega_k,\quad 
\wh\omega_{\{k,l\}}=\frac{1}{b_k-b_l}(\wh\omega_k-\wh\omega_l),\text{ etc.}
\]
We shall use the following elementary lemma in below.
\begin{lem}[Euler]\label{euler.lem}
Let $m\geq1$ be an integer and put $I=\{1,2,\ldots,m\}$.
Let $x_i$ be indeterminates.
Then
\[
\sum_{k\in I}\left(\prod_{j\in I\setminus\{k\}}\frac{1}{x_k-x_j}\right)x_k^r
=\begin{cases}
1&r=m-1\\
0&0\leq r<m-1.
\end{cases}
\]
If $r\geq m$, then
the left hand side is a homogeneous symmetric polynomial of degree $r-m+1$.
\end{lem}
\begin{proof}
Let $F(T)=(T-x_1)\cdots(T-x_m)$ with $T$ another indeterminate. 
Let
\[
\sum_{k=1}^m\frac{1}{F'(x_k)}\frac{1}{T-x_k} =\frac{1}{F(T)}
\]
be the partial fraction decomposition
in $\mathbb{Q}(x_1,\ldots, x_m, T)$. By expanding it in $\mathbb{Q}(x_1, \ldots, x_m)((T^{-1}))$, we have
\[
\sum_{r= 0}^\infty\sum_{k=1}^m\frac{1}{F'(x_k)}x_k^rT^{-r} = T^{-1+m}
\prod_{k=1}^m\left(\sum_{r_k= 0}^\infty x_k^{r_k}T^{-r_k}\right).
\]
Comparing the coefficients of $T^{-r}$, the assertion follows.
\end{proof}


\begin{lem}\label{base.lem1}
Let
$
\wh\omega_I=
G_{1,I}(z)\omega+G_{2,I}(z)D\omega+\cdots+G_{n,I}(z)D^{n-1}\omega.
$
Then,
for any $k\ne l\in I$,
every coefficients of $G_{i,I}(z)$ have no factor $(b_k-b_l)$ in their denominators
(here we think $b_j$'s of being indeterminates).
\end{lem}
\begin{proof}
The factor $(b_k-b_l)$ may appear in the denominators of the term
\[
\left(\prod_{i\in I\setminus\{k\}}\frac{1}{b_k-b_i}\right)\wh\omega_k+
\left(\prod_{i\in I\setminus\{l\}}\frac{1}{b_l-b_i}\right)\wh\omega_l.
\]
This is equal to
\[
\left(\prod_{i\in I\setminus\{k\}}\frac{1}{b_k-b_i}\right)(\wh\omega_k-\wh\omega_l)+
\left(\prod_{i\in I\setminus\{l\}}\frac{1}{b_l-b_i}
+\prod_{i\in I\setminus\{k\}}\frac{1}{b_k-b_i}\right)\wh\omega_l,
\]
and the factor $(b_k-b_l)$
does not appear in the denominators in the first term.
Therefore
 it is enough to see
 \[\prod_{i\in I\setminus\{l\}}\frac{1}{b_l-b_i}
+\prod_{i\in I\setminus\{k\}}\frac{1}{b_k-b_i}
=\frac{1}{b_l-b_k}\left(
\prod_{i\in I\setminus\{k,l\}}\frac{1}{b_l-b_i}
-\prod_{i\in I\setminus\{k,l\}}\frac{1}{b_k-b_i}\right)
\]
has no factor $(b_k-b_l)$ in the denominator.
However, this is easy to see.
\end{proof}

\begin{lem}\label{base.lem2}
For a permutation $\tau$ on $\{1,2,\ldots,n\}$, we define
$\ul b^\tau=(b_{\tau(1)},\ldots,b_{\tau(n)})$ and
\[
\tau(\wh\omega_I)=\tau(\wh\omega_I(\ul a,\ul b))=\wh\omega_{I}(\ul a,\ul b^\tau).
\]
Then we have
\[
\tau(\wh\omega_k)=\wh\omega_{\tau(k)}(\ul a,\ul b),\quad
\tau(\wh\omega_I)=\wh\omega_{\tau(I)}(\ul a,\ul b)
\]
where $\tau(I)=\{\tau(i_1),\ldots,\tau(i_m)\}$.
\end{lem}
\begin{proof}
This follows from the construction of $\wh\omega_I$ and the fact
that 
$F(\ul a;\ul b;z)=F(\ul a;\ul b^\tau;z)$
and $P(\ul a;\ul b)=P(\ul a;\ul b^\tau)$.
\end{proof}

Suppose that $0<b_i<1$ and $b_i\ne b_j$ for any $i\ne j$.
Fix  an integer $s$ such that $1\leq s\leq n$.
We consider the limit $b_1,b_2,\ldots,b_s\to1$, which we simply write $b\to1$.
Put $\ul b_0:=(1,\ldots,1,b_{s+1},\ldots,b_n)$ and 
\[
P_0:=P(\ul a;\ul b_0)=\lim_{b\to 1}P=D^s\prod_{j=s+1}^n(D+b_j-1)-z\prod_{i=1}^n(D+a_i).
\]
The limit $\lim_{b\to1}\wh\omega_I\in \wh H_{P_0}$ 
exists by Lemma \ref{base.lem1}.
For any permutation $\tau$ such that $\tau(i)=i$ for all $i>s$,
\[
\lim_{b\to1}\tau(\wh\omega_I)=\lim_{b\to1}\wh\omega_I\in \wh H_{P_0}
\]
by definition, so that
we have
\begin{equation}\label{omega.perm}
\lim_{b\to1}\wh\omega_{\tau(I)}=\lim_{b\to1}\wh\omega_I\in \wh H_{P_0}
\end{equation}
by Lemma \ref{base.lem2}.

\begin{defn}\label{omega-basis-defn}
Let $0\leq m\leq s$.
We define
\begin{equation}\label{omega-basis}
\wh\omega(m):=\lim_{b\to1}\wh\omega_I\in \wh H_{P_0}
\end{equation}
for a subset $I\subset\{1,2,\ldots,s\}$ such that
$\sharp I=m$. Note $\wh\omega(0)=0$.
This does not depend on either 
the choice of $(b_1,\ldots,b_s)$ that converges to $(1,\ldots,1)$ or
the choice of $I$ by \eqref{omega.perm},
so that it depends only on $(\ul a,\ul b_0)$.
\end{defn}

\begin{lem}\label{base.lem3}
Let $I\subset \{1,2,\ldots,s\}$ be a non-empty subset.
Then $D\wh\omega_I=-\wh\omega_{I\setminus\{m\}}-(b_m-1)\wh\omega_I$ for any $m\in I$.
In particular,
$D\wh\omega(m)=-\wh\omega(m-1)$ and $D^m\wh\omega(m)=0$ for any $m\in\{1,2,\ldots, s\}$.
\end{lem}

\begin{proof}
By the definition of $\wh\omega_I$ and \eqref{wh.omega.lem.eq1}, we have
\[
D\wh\omega_I + (b_m - 1)\wh\omega_I = \sum_{k\in I}\left(\prod_{i\in I\setminus\{k\}}\frac{1}{b_k - b_i}\right)\left((1-b_k)+(b_m-1)\right)\wh\omega_k.
\]
Since, for each $k\in I\setminus\{m\}$,
\[
\prod_{i\in I\setminus\{k\}}\frac{1}{b_k - b_i}\cdot(b_m-b_k) = -\prod_{i\in I\setminus\{k, m\}}\frac{1}{b_k - b_i}
\]
the right hand side is equal to $-\wh\omega_{I\setminus\{m\}}$.
\end{proof}


Next, we shall show that
$\wh\omega(1),\ldots,\wh\omega(s),\wh\omega_{s+1}(\ul a,\ul b_0),\ldots
,\wh\omega_n(\ul a,\ul b_0)$
forms a basis of $\wh H_{P_0}$.
Put
\begin{equation}\label{Deligne-ext}
\wh H_{P_0,\Q_p[[z]]}:=\sum_{i=0}^{n-1} \Q_p[[z]]\,D^i\omega\subset \wh H_{P_0}
\end{equation}
a free $\Q_p[[z]]$-module of rank $n$.
Obviously $\Q_p((z))\ot\wh H_{P_0,\Q_p[[z]]}=\wh H_{P_0}$.
If $0<b_j\leq 1$, then this is the canonical extension 
in the sense of Deligne \cite[5.1]{deligne-book}.
The submodule
\[
\sum_{m=1}^s \Q_p[[z]]\,\wh\omega(m)
+\sum_{j>s} \Q_p[[z]]\,\wh\omega_j(\ul a,\ul b_0)
\subset \wh H_{P_0,\Q_p[[z]]}
\]
is stable under $D$ by Lemma \ref{wh.omega.lem} \eqref{wh.omega.lem.eq1} and
 Lemma \ref{base.lem3}.
\begin{lem}\label{base.lem4}
For $u\in \wh H_{P_0,\Q_p[[z]]}$, let $u|_{z=0}
\in \wh H_{P_0,\Q_p[[z]]}/z \wh H_{P_0,\Q_p[[z]]}$ denote the reduction modulo $z$.
Then 
\[
\wh\omega(m)|_{z=0}=(-1)^{m+1}\left(D^{s-m}\prod_{i>s}(D+b_i-1)\right)\omega|_{z=0},
\quad m=1,2,\ldots,s.
\]
\end{lem}
\begin{proof}
It is enough to show the lemma for $\wh\omega(s)$ by Lemma \ref{base.lem3}.
Let $J_k=\{1,2,\ldots,k\}$.
Recall the definition of $\wh\omega(m)$,
\[
\wh\omega(m)=\lim_{b\to1}\sum_{k=1}^m\left(\prod_{i\in J_m\setminus\{k\}}\frac{1}{b_k-b_i}\right)\wh\omega_k
\]
where 
\[
\wh\omega_k:=y_{n-1}^{[k]}(D+b_k-1)^{n-1}\omega+\cdots+y^{[k]}_1(D+b_k-1)\omega+y_0^{[k]}\omega,
\]
\[
y_{n-1}^{[k]}=(1-z)\check G_k(z),\quad 
y_i^{[k]}=q^{[k]}_{i+1}y^{[k]}_{n-1}-D\bullet y^{[k]}_{i+1},\quad (0\leq i\leq n-2),
\]
\[
\prod_{j=1}^{n}(D+b_j-b_k)-z\prod_{i=1}^{n}(D+1+a_i-b_k)
=(1-z)(D^{n}+q^{[k]}_{n-1}D^{n-1}+\cdots+q^{[k]}_1D+q^{[k]}_0).
\]
Put $\ve_k=b_k-1$ and
$\Delta(z)=(z-\ve_1)\cdots(z-\ve_s)$,
$F(z)=\Delta(z)\prod_{j>s}(z-\ve_j)$.
Then
\[
y_{n-1}^{[k]}(0)=1,\quad y_i^{[k]}(0)=q_{i+1}^{[k]}(0)
=\frac{(-1)^{i+1+n}}{(i+1)!}\frac{d^{i+1}F}{dz^{i+1}}(\ve_k).
\]
We thus have
\[
\wh\omega_k|_{z=0}=\sum_{i=0}^{n-1} \frac{(-1)^{i+1+n}}{(i+1)!}
\frac{d^{i+1}F}{dz^{i+1}}(\ve_k) \cdot(D+\ve_k)^i
\omega|_{z=0}
\]
and hence
\[
\wh\omega(s)|_{z=0}=\lim_{\ve_1,\ldots,\ve_s\to0}
\sum_{k=1}^s\frac{1}{\Delta'(\ve_k)}\sum_{i=0}^{n-1} \frac{(-1)^{i+1+n}}{(i+1)!}
\frac{d^{i+1}F}{dz^{i+1}}(\ve_k) \cdot
(D+\ve_k)^i
\omega|_{z=0}.
\]
Therefore, if we show an equality
\begin{equation}\label{base.lem4.eq1}
\sum_{k=1}^s\frac{1}{\Delta'(\ve_k)}
\sum_{i=0}^{n-1} \frac{(-1)^{i+1+n}}{(i+1)!}(x+\ve_k)^i\frac{d^{i+1}F}{dz^{i+1}}(\ve_k) 
=(-1)^{s+1}\prod_{i>s}(x+\ve_i),
\end{equation}
where $x$ is an indeterminate, this finishes the proof.
However, by the Taylor expansion, we have 
\begin{align*}
\text{left hand side of \eqref{base.lem4.eq1}}
&= \sum_{k=1}^s\frac{1}{\Delta'(\ve_k)}(-1)^n\frac{F(-x)-F(\ve_k)}{x+\ve_k} \\
&= \sum_{k=1}^s\frac{(-1)^n}{\Delta'(\ve_k)}\frac{F(-x)}{x+\ve_k} \\
&= (-1)^{n+1}F(-x)\sum_{k=1}^s\frac{1}{\Delta'(\ve_k)(-x-\ve_k)} \\
&= (-1)^{n+1}F(-x)\frac{1}{\Delta(-x)},
\end{align*}
which implies the desired equality.

\end{proof}

\begin{prop}\label{base.thm}
Suppose that $\{b_j\}_{j>s}$ are pairwise distinct and $b_j\ne1$ for all $j>s$.
Then 
\[\wh\omega(1),\ldots,\wh\omega(s),
\wh\omega_{s+1}(\ul a,\ul b_0),\ldots,\wh\omega_n(\ul a,\ul b_0)
\]
forms a basis of $\wh H_{P_0,\Q_p[[z]]}$. 
\end{prop}
\begin{proof}
By definition, 
$\wh H_{P_0,\Q_p[[z]]}$ is a free $\Q_p[[z]]$-module of rank $n$.
Put $M:=\wh H_{P_0,\Q_p[[z]]}/z\wh H_{P_0,\Q_p[[z]]}$.
We show that
\[
\wh\omega(1)|_{z=0},\ldots,\wh\omega(s)|_{z=0},\,
\wh\omega_{s+1}(\ul a,\ul b_0)|_{z=0},\ldots, \wh\omega_{n}(\ul a,\ul b_0)|_{z=0}\in M
\]
forms a basis of $M$.
Then the assertion follows by Nakayama's lemma.
Note that $M$ is generated by
$\omega|_{z=0}$ as $\Q_p[D]$-module, and it is annihilated by
$D^sQ$ where $Q:=\prod_{j>s}(D+b_j-1)$.
Therefore, one has an isomorphism $M\cong \Q_p[D]/(D^sQ)$, and this
gives rise to an exact sequence
\[
\xymatrix{
0\ar[r]&\bigoplus_{j>s}\ker(D+b_j-1)\ar[r]&M\ar[r]&\Image(Q)\ar[r]&0}.
\]
It follows from Lemmas \ref{base.lem3} and \ref{base.lem4} that 
$\wh\omega(1)|_{z=0},\ldots,\wh\omega(s)|_{z=0}$ span $\Image(Q)$.
Let $j>s$.
The kernel $\ker(D+b_j-1)$ is 1-dimensional, since
$\{b_j\}_{j>s}$ are pairwise distinct and $b_j\ne1$.
One easily shows $\wh\omega_j(\ul a,\ul b_0)|_{z=0}\ne0$ in $M$,
and hence it is a basis of $\ker(D+b_j-1)$
(Lemma \ref{wh.omega.lem}).
Hence $M$ is generated by $\{\wh\omega(m)\}_{1\leq m\leq s}$ and $\{\wh\omega_j(\ul a,\ul b_0)|_{z=0}\}_{j>s}$.
This completes the proof.
\end{proof}

\subsection{Generalized Kedlaya's residue formula}
\begin{defn}\label{Psi-defn}
\begin{itemize}
\item[\rm(1)]
Let $\ul a=(a_1,\ldots,a_n)\in\Z_p^n$ and
$\ul b=(b_1,\ldots,b_n)\in\Z_p^n$ and
put $\psi^{(r)}_p(\ul a,\ul b):=\sum_{i=1}^n\wt\psi^{(r)}_p(a_i)-\wt\psi^{(r)}_p(b_i)$.
We define $\Psi_m(\ul a,\ul b)$ by the following series expansion
\[
\exp\left(\sum_{r=1}^\infty \psi^{(r-1)}_p(\ul a,\ul b)\frac{x^r}r\right)
=\sum_{m=0}^\infty \Psi_m(\ul a,\ul b)x^m.
\] 
\item[\rm(2)]
Let $s\in\{1,2,\ldots,n\}$.
Suppose that $b_j\ne1$ and $b_j^{(1)}\ne 1$ for all $j>s$.
Put
$ \beta_j=1/(b_j-1)$ and 
$\beta^\FS_j=p^{-1}/(b_j^{(1)}-1)$ for $j>s$.
We define $\Psi_m^{[s]}(\ul a,\ul b)$ by
\[
\left(\prod_{j>s}\frac{1-x \beta^\FS_j}{1-x \beta_j}\right)
\exp\left(\sum_{r=1}^\infty \psi^{(r-1)}_p(\ul a,\ul b)\frac{x^r}r\right)
=\sum_{m=0}^\infty \Psi_m^{[s]}(\ul a,\ul b)x^m.
\]
Note that $\Psi_m^{[n]}(\ul a,\ul b)=\Psi_m(\ul a,\ul b)$.
\end{itemize}
\end{defn}

The following is our generalization of Kedlaya's residue formula, which plays a key role in later sections.
\begin{thm}\label{main.formula1}
Let $n\geq s\geq 1$ be integers.
Let $\ul a=(a_1,\ldots,a_n)\in\Z_{(p)}^n$ and 
$\ul b_0=(1,\ldots,1,b_{s+1},\ldots,b_n)\in\Z_{(p)}^n$ satisfy that
$0<a_i,b_j<1$ and
$a_i-b_j\not\in\Z$ for every $i$ and $j>s$, and 
that $b_{s+1},\ldots,b_n$ are pairwise distinct.
Let $\sigma(z)=z^p$, and let
\[
\xymatrix{
\Phi:H_{P(\ul a^{(1)};\ul b_0^{(1)})}^\dag\ar[r]&H_{P(\ul a;\ul b_0)}^\dag
}\]
be the $\sigma$-linear Frobenius intertwiner, that is unique
up to $\Q_p^\times$.
Let
\[
\{\wh\omega(1),\ldots,\wh\omega(s),\wh\omega_{s+1},\ldots,\wh\omega_n\},\quad
\{\wh\omega^\FS(1),\ldots,\wh\omega^\FS(s),\wh\omega^\FS_{s+1},\ldots,\wh\omega^\FS_n\}
\]
be the bases of 
$\wh H_{P(\ul a;\ul b_0)}$ and $\wh H_{P(\ul a^{(1)};\ul b^{(1)}_0)}$ 
in Proposition \ref{base.thm}
respectively, where we simply write
$\wh\omega_k=\wh\omega_k(\ul a,\ul b_0)$
and 
$\wh\omega^\FS_k=\wh\omega_k(\ul a^{(1)},\ul b_0^{(1)})$.
Put
\[
\gamma_1
:=p^{-1}\prod_{i=1}^{n}\frac{\Gamma_p(b_i)}{\Gamma_p(a_i)},\quad
\gamma'_1
:=p^{-s}\prod_{i>s}\frac{1-b_i^{(1)}}{1-b_i}\prod_{i=1}^{n}\frac{\Gamma_p(b_i)}{\Gamma_p(a_i)},
\]
and $\gamma_m=\gamma_1$ and $\gamma'_m=\gamma'_1$ for $1\leq m\leq s$.
For $k>s$, let $\gamma_k$ be the constants \eqref{Kedlaya-gammak}
for $(\ul a,\ul b_0)$, and 
put
\[
\gamma'_k:=
\gamma_k\left(\frac{b_k^{(1)}-1}{b_k-1}\right)^s\prod_{i>s,\,i\ne k}\frac{b_k^{(1)}-b_i^{(1)}}{b_k-b_i}.
\]
Put $\mu_k=p(1-b^{(1)}_k)-(1-b_k)$ for $k>s$.
Fix the intertwiner $\Phi$ so that $\Phi(\wh\omega^\FS(1))=\wh\omega(1)$.
Then, we have
\begin{equation}\label{main.formula1-eq1}
p^{-m+1}
\Phi(\wh\omega^\FS(m))=\sum_{i=1}^m\Psi^{[s]}_{m-i}(\ul a,\ul b_0)\,\wh\omega(i),\quad
(1\leq m\leq s),
\end{equation}
\begin{equation}\label{main.formula1-eq2}
\Phi(\wh\omega^\FS_k)=(\gamma'_k/\gamma'_1)z^{\mu_k}\wh\omega_k,\quad (k>s).
\end{equation}
\end{thm}
\begin{rem}
In case of $s=n$ (i.e. $\ul b_0=(1,\ldots,1)$), the condition on $(\ul a,\ul b_0)$ in 
Theorem \ref{main.formula1}
is simply that $0<a_i<1$ 
for every $i$.
\end{rem}

\subsection{Proof of Theorem \ref{main.formula1}}\label{proof.main.formula.sect}
We shall derive \eqref{main.formula1-eq1} and \eqref{main.formula1-eq2}
by taking the limit of Kedlaya's
formula (Theorem \ref{Kedlaya}).
The formula \eqref{main.formula1-eq2} is almost straightforward, 
while the proof of \eqref{main.formula1-eq1} is non-trivial.

\medskip

Let $\ul a=(a_1,\ldots,a_n)$ and
$\ul b=(b_1,\ldots,b_n)$ be parameters satisfying \textbf{Hypothesis} \ref{Hypothesis}.
Let $s$ be an integer such that $1\leq s\leq n$,
and we suppose $b_1\equiv\cdots\equiv b_s\equiv 1$ mod $p$.
Let
$P(\ul a;\ul b)$ and 
$P(\ul a^{(1)};\ul b^{(1)})$ 
be the hypergeometric differential operators as in \eqref{HG.diff.op}.
Let 
$\{\wh\omega_{k,b}\}_{1\leq k\leq n}$
and 
$\{\wh\omega^\FS_{k,b}\}_{1\leq k\leq n}$
be the elements in $\wh H_{P(\ul a;\ul b)}$ and $\wh H_{P(\ul a^{(1)};\ul b^{(1)})}$
respectively defined in \eqref{wh.omega.defn}.
Let
\begin{equation}\label{P0-eq1}
\xymatrix{
\Phi_b:\wh H_{P(\ul a^{(1)};\ul b^{(1)})}\ar[r]&\wh H_{P(\ul a;\ul b)},}
\quad\xymatrix{
\Phi:\wh H_{P(\ul a^{(1)};\ul b^{(1)}_0)}\ar[r]&\wh H_{P(\ul a;\ul b_0)}}
\end{equation}
be the $\sigma$-linear Frobenius intertwiners, which are 
uniquely determined up to $\Q_p^\times$.
To remove the ambiguity of $\Q_p^\times$, we fix $\Phi_b$ as follows, and
take $\Phi$ to be the
limit of $\Phi_b$ as $b_1,b_2,\ldots,b_s\to 1$.
Let
\[
c_{k,b}:=\prod_{1\leq i\leq n,\,i\ne k}(b_k-b_i) 
,\quad c^\FS_{k,b}:=\prod_{1\leq i\leq n,\,i\ne k}(b^{(1)}_k-b^{(1)}_i)
\]
and let $\gamma_{1,b},\ldots,\gamma_{n,b}$ be
the constants \eqref{Kedlaya-gammak} for $(\ul a,\ul b)$.
Put $\gamma'_{k,b}:=c_{k,b}^{-1}c^\FS_{k,b}\gamma_{k,b}$.
Then, we choose $\Phi_b$ to be the Frobenius intertwiner
satisfying 
\begin{equation}\label{P0-eq2}
\Phi_b(\wh\omega^\FS_{k,b})=\gamma'_{k,b}z^{\mu_k}\,
\wh\omega_{k,b},\quad (k=1,2,\ldots,n).
\end{equation}
Take the limit $b_1,\ldots,b_s\to1$, then one sees $\gamma'_{k,b}\to\gamma'_k$
for $k=1,2,\ldots,n$.
Indeed, if $b_k\equiv1$ mod $p$, then this immediately follows from
Lemma \ref{gamma.formula}.
If $b_k\not\equiv1$ mod $p$, then, using the notation in the proof of 
Lemma \ref{gamma.formula}, one can show that
$A(x)=b_k^{(1)}-x^{(1)}$ for $x\in \Z_{(p)}\cap (0,1)$ such that $x\equiv 1$ mod $p$, so that
one has
\[
\gamma_{k,b}=\prod_{i=1}^nA(a_i)\cdot\prod_{i=1}^s(b_k^{(1)}-b_i^{(1)})^{-1}\cdot
\prod_{j>s}A(b_j)^{-1}\cdot 
\prod_{i=1}^n\frac{\Gamma_p(b_i-b_k+1)}{\Gamma_p(a_i-b_k+1)}\to \gamma_k.
\]
We thus have
\[
\Phi(\wh\omega^\FS_k)=\gamma'_kz^{\mu_k}\,
\wh\omega_k,\quad (k=1,2,\ldots,n),
\]
and hence if we let the $\Phi$ in 
Theorem \ref{main.formula1} to be $(\gamma'_1)^{-1}\Phi$, then
\eqref{main.formula1-eq2} holds.

\bigskip

Next, we show \eqref{main.formula1-eq1}.
\begin{lem}\label{gamma.Psi.formula1}
Let $m$ be an integer such that $1\leq m\leq s$.
Put $J_m=\{1,2,\ldots,m\}$.
Then
\begin{equation}\label{gamma.Psi1.eq1}
\lim_{b_1,\ldots,b_s\to1}
\sum_{k=1}^m\left(\prod_{j\in J_m\setminus\{k\}}\frac{1}{b_k-b_j}\right)
\frac{\gamma_{k,b}}{\gamma_{1,b}}
=\Psi_{m-1}(\ul a,\ul b_0).
\end{equation}
\end{lem}
\begin{proof}
Put $\ve_k=b_k-1$ for $k=1,2,\ldots,s$.
It follows from Lemma \ref{gamma.formula.cor} and Theorem \ref{log.beta} that
one has
\begin{align*}
\frac{\gamma_{k,b}}{\gamma_{1,b}}
&=
\prod_{i=1}^n\frac{B_p(b_i,-\ve_1)}{B_p(b_i,-\ve_k)}\frac{B_p(a_i,-\ve_k)}{B_p(a_i,-\ve_1)}\\
&=
\exp\left(\sum_{r=1}^\infty
\left(\sum_{i=1}^n\wt\psi^{(r-1)}_p(a_i)-\wt\psi^{(r-1)}_p(b_i)\right)\frac{\ve_k^r}{r}
-\left(\sum_{i=1}^n\wt\psi^{(r-1)}_p(a_i)-\wt\psi^{(r-1)}_p(b_i)\right)
\frac{\ve_1^r}{r}
\right)\\
&=
\exp\left(\sum_{r=1}^\infty
\psi^{(r-1)}_p(\ul a,\ul b)
\left(\frac{\ve_k^r}{r}-\frac{\ve_1^r}{r}
\right)\right)
\end{align*}
Let us put
\[
F(\ve_1,\ldots,\ve_s,x):=
\exp\left(\sum_{r=1}^\infty
(\sum_{i=1}^s\wt\psi^{(r-1)}_p(a_i)-\wt\psi^{(r-1)}_p(b_i))\left(\frac{x^r}{r}-\frac{\ve_1^r}{r}
\right)\right)
\]
an analytic function of variables $\ve_1,\ldots,\ve_s,x$ where
we think $a_i$ and $b_j$ ($j>s$) of being constants.
This has a series expansion that
converges
on $\{|\ve_i|\ll1,|x|\ll1\}_i$ (Lemma \ref{polygamma}).
Put
\[
G(\ve_1,\ldots,\ve_s;x_1,\ldots,x_m):=
\sum_{k=1}^m\left(\prod_{j\in J_m\setminus\{k\}}\frac{1}{x_k-x_j}\right)
F(\ve_1,\ldots,\ve_s,x_k).
\]
By Lemma \ref{euler.lem}, this 
is also an analytic function that has a series expansion
converging
on $\{|\ve_i|\ll1,|x_k|\ll1\}_{i,k}$.
We have
\[
\text{left hand side of \eqref{gamma.Psi1.eq1}}=\lim_{\ve_1,\ldots,\ve_s\to0}
G(\ve_1,\ldots,\ve_s;\ve_1,\ldots,\ve_m).
\]
However, this is equal to
\begin{align*}
&
\lim_{\ve_1,\ldots,\ve_s\to0}
\lim_{x_1,\ldots,x_m\to0}
G(\ve_1,\ldots,\ve_s;x_1,\ldots,x_m)\\
=&
\lim_{x_1,\ldots,x_m\to0}
G(0,\ldots,0;x_1,\ldots,x_m)\\
=&
\lim_{x_1,\ldots,x_m\to0}
\sum_{k=1}^m\left(\prod_{j\in J_m\setminus\{k\}}\frac{1}{x_k-x_j}\right)
\exp\left(\sum_{r=1}^\infty
\psi^{(r-1)}_p(\ul a,\ul b_0)
\frac{x_k^r}{r}
\right)\\
=&
\lim_{x_1,\ldots,x_m\to0}
\sum_{k=1}^m\left(\prod_{j\in J_m\setminus\{k\}}\frac{1}{x_k-x_j}\right)
\left(\sum_{r=1}^\infty
\Psi_{r-1}(\ul a,\ul b_0)x_k^r
\right)\\
=&\Psi_{m-1}(\ul a,\ul b_0),
\end{align*}
where the last equality follows from Lemma \ref{euler.lem}.
This completes the proof.
\end{proof}

\begin{lem}\label{gamma.Psi.formula2}
For $m=1,2,\ldots,s$, we have
\[
\lim_{b_1,\ldots,b_s\to1}
\sum_{k=1}^m\left(\prod_{i\in J_m\setminus\{k\}}\frac{1}{b_k-b_i}\right)
\frac{\gamma'_{k,b}}{\gamma'_{1,b}}
=\Psi_{m-1}^{[s]}(\ul a,\ul b_0).
\]
\end{lem}
\begin{proof}
Put $\ve_k=b_k-1$ for $1\leq k\leq s$. 
We have
\[
\frac{c^\FS_{k,b}}{c_{k,b}}=p^{-s+1}\frac{(b_k^{(1)}-b_{s+1}^{(1)})\cdots(b_k^{(1)}-b_n^{(1)})}
{(b_k-b_{s+1})\cdots(b_k-b_n)}=p^{-s+1}\frac{(\ve_k/p+1-b_{s+1}^{(1)})\cdots(\ve_k/p+1-b_n^{(1)})}
{(\ve_k+1-b_{s+1})\cdots(\ve_k+1-b_n)},
\]
and
\[
\frac{\gamma'_{k,b}}{\gamma'_{1,b}}=\frac{\gamma_{k,b}}{\gamma_{1,b}}
\frac{c^\FS_{k,b}}{c_{k,b}}\frac{c_{1,b}}{c^\FS_{1,b}}=\frac{\gamma_{k,b}}{\gamma_{1,b}}
\left(\prod_{i>s}\frac{1-\ve_k \beta^\FS_i}{1-\ve_k \beta_i}\right)
\left(\prod_{i>s}\frac{1-\ve_1 \beta^\FS_i}{1-\ve_1 \beta_i}\right)^{-1}.
\]
The rest goes in a similar way to the proof of Lemma \ref{gamma.Psi.formula1}.
\end{proof}

\begin{lem}\label{Frob.gamma.lem}
Put $J_m=\{1,2,\ldots,m\}$ and
\[\gamma'(m):=\lim_{b_1,\ldots,b_s\to1}
\sum_{k=1}^m\left(\prod_{i\in J_m\setminus\{k\}}\frac{1}{b_k-b_i}\right)\gamma'_{k,b}.
\]
Note $\gamma'(1)=\gamma'_1$.
Then
\[
p^{-m+1}\Phi(\wh\omega^\FS(m))
=\gamma'(1)\wh\omega(m)+\gamma'(2)\wh\omega(m-1)
+\cdots+\gamma'(m)\wh\omega(1)
\]
for $m=1,2,\ldots,s$.
\end{lem}
\begin{proof}
Put
\[
A_b=p^{-m+1}\Phi_b\left(\sum_{k=1}^m\left(
\prod_{i\in J_m\setminus\{k\}}\frac{1}{b^{(1)}_k-b^{(1)}_i}\right)\wh\omega_{k,b}^\FS\right).
\]
Since $b^{(1)}_k-b^{(1)}_i=p(b_k-b_i)$ for $i,k\in J_s$, 
we have
\[
A_b=\sum_{k=1}^m\left(\prod_{i\in J_m\setminus\{k\}}\frac{1}{b_k-b_i}\right)
\Phi_b(\wh\omega^\FS_{k,b})=
\sum_{k=1}^m\left(\prod_{i\in J_m\setminus\{k\}}\frac{1}{b_k-b_i}\right)
\gamma'_{k,b}\wh\omega_{k,b}.
\]
On the other hand, we have
\begin{align*}
 &\sum_{\alpha=1}^m\left[\sum_{1\le k\le \alpha}\left(\prod_{1\le i\le \alpha, i\neq k} \frac{1}{b_k - b_i}\right)\gamma'_{k, b}\right]\left[\sum_{\alpha\le l\le m}\left(\prod_{\alpha\le j\le m, j\neq l} \frac{1}{b_l - b_j}\right)\wh\omega_{l, b}\right]\\
=&\sum_{1\le k\le l\le m}\left[\sum_{k\le \alpha\le l}\left(\prod_{1\le i\le \alpha, i\neq k} \frac{1}{b_k - b_i}\right)\left(\prod_{\alpha\le j\le m, j\neq l} \frac{1}{b_l - b_j}\right)\right]\gamma'_{k, b}\wh\omega_{l, b}\\
=& \sum_{1\le k\le m}\left(\prod_{1\le i\le m, i\neq k} \frac{1}{b_k - b_i}\right)\gamma'_{k, b}\wh\omega_{k, b} = A_b,
\end{align*}
where the last equality follows from Lemma \ref{lemma 4.14} below. 
Taking the limit $b_1,\ldots,b_s\to1$, we have the lemma.
\end{proof}

\begin{lem}\label{lemma 4.14}
For $1\le k\le l \le m$, we have
\[
\sum_{k\le \alpha\le l}\left(\prod_{1\le i\le \alpha, i\neq k} \frac{1}{b_k - b_i}\right)\left(\prod_{\alpha\le j\le m, j\neq l} \frac{1}{b_l - b_j}\right)
 = \prod_{1\le i\le m, i\neq k} \frac{1}{b_k - b_i}\delta_{k, l},
\]
where $\delta_{k, l}$ denotes Kronecker's delta.
\end{lem}
\begin{proof}
The case where $l-k\le 1$ can be proven by direct calculation. If $l>k+1$, we can show the vanishing of the left hand side by replacing
\[
\left(\prod_{1\le i\le \alpha, i\neq k} \frac{1}{b_k - b_i}\right)\left(\prod_{\alpha\le j\le m, j\neq l} \frac{1}{b_l - b_j}\right)
\]
with
\[
\left(\prod_{1\le i < \alpha, i\neq k} \frac{1}{b_k - b_i}\right)\cdot\Biggl(\frac{1}{b_k - b_l}\frac{1}{b_l-b_\alpha}-\frac{1}{b_k - b_\alpha}\frac{1}{b_k - b_l}\Biggr)\cdot\left(\prod_{\alpha < j\le m, j\neq l} \frac{1}{b_l - b_j}\right)
\]
for $k<\alpha<l$.
\end{proof}

Let the $\Phi$ in 
Theorem \ref{main.formula1} be $(\gamma'_1)^{-1}\Phi$ as before. Then
\eqref{main.formula1-eq1}
follows from \eqref{P0-eq2} and Lemmas \ref{gamma.Psi.formula2}
and \ref{Frob.gamma.lem}.
This completes the proof of Theorem \ref{main.formula1}.

\subsection{Scalar extension and Changing the $p$-th Frobenius}\label{changing.sect}
Let $W=W(k)$ be the Witt ring of a perfect field $k$ of characteristic $p>0$,
and $K=\Frac(W)$ the fractional field.
The $\sigma$-linear Frobenius intertwiner $\Phi=\Phi_\sigma$ on $H_{P(\ul a;\ul b)}$
naturally extends on the scalar extension
$H_{P(\ul a;\ul b),K}:=K\ot_{\Q_p}H_{P(\ul a;\ul b)}$, which we write by the same notation.
Let $\tau$ be another $p$-th Frobenius on $K[z,(z-z^2)^{-1}]^\dag$.
One can associate the $\tau$-linear Frobenius intertwiner $\Phi_\tau$
to $\Phi$ in 
the canonical way by the theory of $F$-isocrystals, that is explicitly given 
as follows (\cite[6.1]{EK}, 
\cite[17.3.1]{Ke-DE}),
\begin{equation}\label{change-Frob}
\Phi_\tau(v)-\Phi_\sigma(v)=
\sum_{i=1}^\infty\frac{(z^\tau-z^\sigma)^i}{i!}
\Phi_\sigma\partial^i v,\quad \partial:=\frac{d}{dz}.
\end{equation}
We consider a $p$-th Frobenius $\tau$ on $K[z,(z-z^2)^{-1}]^\dag$
given by $\tau(z)=cz^p$ for some $c\in 1+pW$.
We have
\begin{equation}\label{changing.eq1}
\frac{(z^\tau-z^\sigma)^i}{i!}
\Phi_\sigma\partial^i v=
\frac{(c-1)^iz^{ip}}{i!}
\Phi_\sigma(z^{-i}\cdot z^i\partial^i v)
=
\frac{(c-1)^i}{i!}
\Phi_\sigma(D(D-1)\cdots(D-i+1) v).
\end{equation}
Let $v=\wh\omega^\FS_k$ and apply \eqref{changing.eq1}.
Since $D\wh\omega^\FS_k=(1-b^{(1)}_k)\wh\omega^\FS_k$ 
by Lemma \ref{wh.omega.lem},
we have
\[
\Phi_\tau(\wh\omega^\FS_k)=
\sum_{i=0}^\infty(c-1)^i
\binom{1-b^{(1)}_k}{i}\Phi_\sigma(\wh\omega^\FS_k)
=c^{1-b^{(1)}_k}\,\Phi_\sigma(\wh\omega^\FS_k).
\]
Let $\log:1+p\Z_p\to p\Z_p$ be the logarithmic function
defined by $\log(1-px)=-\sum_{n\geq 1}p^nx^n/n$.
Let $v=\wh\omega^\FS(m)$ and apply \eqref{changing.eq1}.
Then we have
\begin{align*}
\Phi_\tau(\wh\omega^\FS(m))
&=\sum_{i=0}^\infty\frac{(c-1)^i}{i!}\Phi_\sigma\left(D(D-1)\cdots(D-i+1)\wh\omega^\FS(m)\right)\\
&=\Phi_\sigma\left(
\sum_{i=0}^\infty\frac{(F^{-1}(c)-1)^i}{i!}D(D-1)\cdots(D-i+1)\wh\omega^\FS(m)\right)\\
&=\Phi_\sigma\left(\sum_{i=0}^\infty(\log F^{-1}(c))^i\frac{D^i}{i!}\wh\omega^\FS(m)\right)\\
&=\Phi_\sigma\left(\sum_{i=0}^{m-1}\frac{(-\log F^{-1}(c))^i}{i!}\wh\omega^\FS(m-i)\right)&(\text{Lemma \ref{base.lem3})}\\
&=\sum_{i=0}^{m-1}\frac{(-\log c)^i}{i!}\Phi_\sigma(\wh\omega^\FS(m-i)).
\end{align*}
Here, the third equality follows from the identity in $\Q[[X, Y]]$
\[
\sum_{i=0}^\infty(\log(1+Y))^i\frac{X^i}{i!}  = \sum_{i=0}^\infty \binom{X}{i}Y^i
\]
(Note that it is a formal expression of the identity ``$\exp(X\log (1+Y)) = (1+Y)^X$'').
Combining the above with Theorem \ref{main.formula1}, we have the following.
\begin{thm}\label{main-c}
Let the notation and setting be as in Theorem \ref{main.formula1}.
Let $\tau(z)=cz^p$ with $c\in 1+pW$, and $\Phi_\tau:H_{P(\ul a^{(1)};\ul b^{(1)}),K}\to
H_{P(\ul a;\ul b),K}$ the 
$\tau$-linear Frobenius intertwiner, 
where $H_{P(\ul a;\ul b),K}:=K\ot_{\Q_p}H_{P(\ul a;\ul b)}$.
Then
\begin{equation}\label{main.c-eq1}
p^{-m+1}\Phi_\tau(\wh\omega^\FS(m))=\sum_{k=1}^m\left(
\sum_{i+j=m-k}\frac{(-p^{-1}\log c)^i}{i!}
\Psi^{[s]}_j(\ul a,\ul b_0)\right)\wh\omega(k)
\end{equation}
for $1\leq m\leq s$, and
\begin{equation}\label{main.c-eq2}
\Phi_\tau(\wh\omega^\FS_k)=c^{1-b^{(1)}_k}(\gamma'_k/\gamma'_1)z^{\mu_k}\wh\omega_k,\quad 
(k>s).
\end{equation}
\end{thm}
\begin{rem}
The coefficients of $\wh\omega(k)$ in \eqref{main.c-eq1} arise
as the coefficients of 
the following series,
\begin{align*}
&\left(\prod_{i>s}\frac{1-x \beta^\FS_i}{1-x \beta_i}\right)
\exp\left(-(p^{-1}\log c)\, x+\sum_{r=1}^\infty \psi^{(r-1)}_p(\ul a,\ul b_0)\frac{x^r}r\right)\\
=&\left(\sum_{i=0}^\infty\frac{(-p^{-1}\log c)^i}{i!}x^i\right)
\left(\sum_{m=0}^\infty \Psi^{[s]}_m(\ul a,\ul b_0)x^m\right)
=\sum_{m=0}^\infty\left(\sum_{i+j=m}\frac{(-p^{-1}\log c)^i}{i!} \Psi^{[s]}_j(\ul a,\ul b_0)\right)x^m.
\end{align*}
\end{rem}
\section{Frobenius on log-crystalline cohomology and $p$-adic $L$-values}
\label{application.sect}
Let $W=W(k)$ be the Witt ring of a perfect field $k$ of characteristic $p>0$,
and $K=\Frac(W)$ the fractional field. 
For a projective smooth morphism $f:X\to S$ of smooth $W$-schemes with $S$ affine, 
let
$H^*_\rig(X_k/S_k)$ denote the {\it rigid cohomology} (we refer to \cite{LS} 
for a general reference of rigid cohomology, but we use the notation in \cite[\S 2.5]{AM}),
where $X_k:=X\times_Wk$ and $S_k:=S\times_Wk$.
The rigid cohomology has the comparison 
\[
H^*_\rig(X_k/S_k)
\cong \O(S_K)^\dag\ot H^*_\dR(X_K/S_K)
\]
with de Rham cohomology (cf. \cite[Proposition 2.6]{AM}), where $X_K=X\times_WK$ etc.
For a $p$-th Frobenius $\tau$ on $\O(S_K)^\dag$,
the $\tau$-linear Frobenius map 
\[\xymatrix{
\Phi_{\rig,\tau}:H^*_\rig(X_k/S_k)\ar[r]& H^*_\rig(X_k/S_k)}\]
is defined in the canonical way, giving a $\tau$-linear Frobenius intertwiner.
If $S=\Spec W[z,(z-z^2)^{-1}]$ and the Picard-Fuchs equation
of the de Rham cohomology is isomorphic to
a hypergeometric equation, then one can obtain an explicit matrix representation
of the
Frobenius $\Phi_{\rig,\tau}$ thanks to Theorems \ref{main.formula1} or \ref{main-c}.
However the Frobenius intertwiner is only determined ``up to multiplication by $K^\times$'',
and there remains issue to remove the ambiguity of
$K^\times$.
In this section, we discuss this issue for
 (i) the Dwork pencils and (ii) families which have
totally degenerating fibers at $z=0$. 
For the latter, we employ the log crystalline cohomology.

\subsection{Frobenius on Katz' generalized Dwork pencils, \cite{Katz}}\label{Katz}
Let $n\geq 2$ be an integer and let
$d> n$ and $w_i>0$ be integers such that $w_0+w_1+\cdots +w_n=d$
and $\gcd(w_0,\ldots,w_n)=1$.
According to Katz \cite{Katz}, we consider a homogeneous equation
 \[
w_0x_0^d+w_1x_1^d+\cdots+w_nx_n^d-\lambda dx_0^{w_0}\cdots x_n^{w_n}
\]
with a parameter $\lambda$. This defines the pencil of hypersurfaces, 
which we call the {\it generalized Dwork pencil}.
The hypersurface is smooth if and only if 
$\lambda^d\ne1$ (\cite[Lemma 2.1]{Katz}).
Choose integers $b_i$ such that $\sum_i b_iw_i=1$, and put $X_i:=\lambda^{b_i}x_i$.
Then, the above equation turns out to be
$
w_0\lambda^{-db_0}X_0^d+\cdots+w_n
\lambda^{-db_n}X_n^d-dX_0^{w_0}\cdots X_n^{w_n}.
$
Furthermore we put $z:=\lambda^{-d}$, then the equation is
\[Q_z=w_0z^{b_0}X_0^d+\cdots+w_n
z^{b_n}X_n^d-dX_0^{w_0}\cdots X_n^{w_n}.
\]
Let $p$ be a prime which is prime to $d$ and $w_0,\ldots,w_n$.
Let $W=W(k)$ be the Witt ring of a perfect field $k$ of characteristic $p$,
and $K=\Frac(W)$ the fractional field. 
Let $f:X\to S:=\Spec W[z,(z-z^2)^{-1}]$ be the projective smooth $W$-morphism
whose general fiber is the hypersurface $Q_z=0$.
Write $X_K=X\times_WK$ and $X_k=X\times_Wk$ etc.
Let $\mu_d$ be the group of $d$-th roots of unity in $\ol K^\times$.
The group $G=\{(\zeta_0,\ldots,\zeta_n)\in \mu_d^{n+1}\mid\zeta_0^{w_0}\cdots\zeta_n^{w_n}=1\}$ acts on $X$ in the following way,
\[
(\zeta_0,\ldots,\zeta_n)(X_0,\ldots,X_n)=(\zeta_0X_0,\ldots,\zeta_nX_n).
\]
Let
\[
 \Omega=\sum_{i=0}^n(-1)^iX_i\,dX_0\wedge\cdots
\wedge\wh{dX_i}\wedge\cdots\wedge dX_n\in \vG(\P^n,\Omega^n_{\P^n/S}(X))
\]
and
\[
\omega_X:=\Res\,\frac{X_0^{w_0-1}\cdots X_n^{w_n-1}}{Q_z}
\Omega\in \vg(X,\Omega^{n-1}_{X/S})
\]
where $\Res$ is the residue map.
Let
\begin{align}
\ul a&=\left(\frac1d,\frac2d,,\ldots,\frac{d-1}{d}\right)\\
\ul b&=\left(\overbrace{1,\ldots,1}^n,
\frac1{w_0},\ldots,\frac{w_0-1}{w_0},\ldots,\frac{1}{w_n},\ldots,\frac{w_n-1}{w_n}\right)\label{b}
\end{align}
be $(d-1)$-tuples,
and let $(\ul a',\ul b')$ be the cancelled lists (namely, if $\ul a$ and $\ul b$ contain 
common entries, then eliminate them).
Then the arrow 
\begin{equation}\label{Katz.eq1}
\xymatrix{
H_{P(\ul a';\ul b')}=
\cD/\cD P(\ul a';\ul b')\ar[r]^-\cong&H^{n-1}_\dR(X_K/S_K)^G,\quad \theta\longmapsto \theta \omega_X}
\end{equation}
is bijective (\cite[Theorem 8.4]{Katz}), where $(-)^G$ denotes the invariant part by
the action of $G$.
Let $c\in 1+pW$ and $\tau$ the $p$-th Frobenius on $\O(S_K)^\dag$ such that
$\tau(z)=cz^p$.
Let $H^*_\rig(X_k/S_k)$ be the rigid cohomology, and 
let 
\[\xymatrix{
\Phi_{\rig,\tau}:H^*_\rig(X_k/S_k)\ar[r]& H^*_\rig(X_k/S_k)}\]
be the $\tau$-linear Frobenius map. 
We have the comparison isomorphism
\begin{equation}\label{Katz.eq2}
\O(S_K)^\dag\ot H^*_\dR(X_K/S_K)\cong 
H^*_\rig(X_k/S_k)
\end{equation}
with de Rham cohomology (cf. \cite[Proposition 2.6]{AM}).
Here we note that $(\ul a')^{(1)}=\ul a'$ and $(\ul b')^{(1)}=\ul b'$ up to permutation, so that
one has $P(\ul a';\ul b')=P((\ul a')^{(1)};(\ul b')^{(1)})$.

\medskip

The following theorem is a generalization of \cite{shapiro} ($d=n+1=5$, $w_i=1$) and
\cite[Theorem 1.4]{BV} ($d=n+1$ arbitrary, $w_i=1$).
\begin{thm}\label{Katz-thm}
Suppose that $\gcd(w_i,w_j)=1$ for every $i\ne j$.
Let $\ul a'$ and $\ul b'$ be the cancelled lists as above, and 
write $\ul a'=(a_1,\ldots,a_{d'})$ and
$\ul b'=(b_1,\ldots,b_{d'})$ where we note $b_1=\cdots=b_s=1$ $(s:=n)$ by \eqref{b}.
Let $\{\wh\omega(1),\ldots,\wh\omega(s),\wh\omega_{s+1},\ldots,\wh\omega_{d'}\}$
be the basis of $\wh H_{P(\ul a';\ul b')}$ as in Proposition \ref{base.thm},
and we think them of being a basis of $H^{n-1}_\rig(X_k/S_k)$
by \eqref{Katz.eq1} and \eqref{Katz.eq2}.
Put $\mu_k=b_k-1-p(b_k^{(1)}-1)$ 
and let $\gamma'_k$ be the constants 
as in Theorem \ref{main.formula1} for $(\ul a,\ul b_0)=(\ul a',\ul b')$.
Then there is a root $\ve\in\Q_p^\times$ of unity
such that
\begin{align}\label{Katz-thm-eq1}
&p^{-m+1}\Phi_{\rig,\tau}(\wh\omega(m))=\ve\sum_{k=1}^m\left(
\sum_{i+j=m-k}\frac{(-p^{-1}\log c)^i}{i!}
\Psi^{[s]}_j(\ul a',\ul b')\right)\wh\omega(k),&(1\leq m\leq s),\\
\label{Katz-thm-eq2}
&\Phi_{\rig,\tau}(\wh\omega_k)=\ve 
c^{1-b^{(1)}_k}(\gamma'_k/\gamma'_1)z^{\mu_k}\wh\omega_k, 
&(k>s).
\end{align}
Moreover, we have $\Psi^{[s]}_m(\ul a',\ul b')=\Psi^{[s]}_m(\ul a,\ul b)$ 
and this lies in the subring
\[
\Q[L, \zeta_p(3), \zeta_p(5),\ldots,\zeta_p(2l+1)]\subset \Q_p
\]
where we put $L:=-dp^{-1}\log(d^{p-1})+\sum_{w_i>1}w_ip^{-1}\log(w_i^{p-1})$ and $l:=\lfloor m/2\rfloor$.
\end{thm}
\begin{proof}
We first note that $X/S$ is defined over $\Z_p$ and every $\wh\omega(m)$ and $\wh\omega_k$
are also defined over $\Q_p$.
Then \eqref{Katz-thm-eq1} and \eqref{Katz-thm-eq2} can be deduced from the case $W=\Z_p$ and
$c=1$ by the same discussion in \S \ref{changing.sect}.
Therefore we may assume that $W=\Z_p$ and $c=1$.

By the assumption that $\gcd(w_i,w_j)=1$ for $i\ne j$,
the entries 
$b_{s+1},\ldots, b_{d'}$ are pairwise distinct, so that
Theorem \ref{main.formula1} is applicable to the Frobenius intertwiner on $H_{P(\ul a';\ul b')}$.
Therefore there is a constant $\ve\in \Q_p^\times$ such that 
\begin{align}
&p^{-m+1}\Phi_{\rig,\tau}(\wh\omega(m))=\ve\sum_{k=1}^m
\Psi_{m-k}^{[s]}(\ul a',\ul b')\wh\omega(k)\label{Katz-thm-eq1-pre}
&(1\leq m\leq s)\\
&\Phi_{\rig,\tau}(\wh\omega_k)=\ve 
(\gamma'_k/\gamma'_1)z^{\mu_k}\wh\omega_k,
&(k>s)\label{Katz-thm-eq2-pre}.
\end{align}
We show that $\ve$ is a root of unity.
Choose $\alpha\in \Z_p^\times$ such that $\alpha^p=\alpha$ and 
$\alpha\not \equiv1$ mod $p$.
Let $X_\alpha$ be the fiber of $X/S$ at $z=\alpha$.
Then $\Phi_{\rig,\tau}$ induces the Frobenius $\Phi_{\rig,\alpha}$ on 
$H^*_\rig(X_{\alpha,\F_p})\cong H^*_\dR(X_\alpha/\Q_p)$.
The cup-product induces a non-degenerate pairing
$H^{n-1}_\dR(X_\alpha/\Q_p)^G\ot H^{n-1}_\dR(X_\alpha/\Q_p)^G\to \Q_p(1-n)$,
and this
is compatible with respect to
the action of $\Phi_{\rig,\alpha}$, where
$\Phi_{\rig,\alpha}$ acts on $\Q_p(j)$ by multiplication by $p^{-j}$.
Therefore one has an isomorphism
$(\det [H^{n-1}_\dR(X_\alpha/\Q_p)^G])^{\ot 2}\cong \Q_p(d'(1-n))$, and this implies that
the product of eigenvalues of $\Phi_{\rig,\alpha}$ on $H^{n-1}_\dR(X_\alpha/\Q_p)^G$
is $\pm p^{\frac12(n-1)d'}$
(note that $d'$ is the rank of $H^{n-1}_\dR(X_\alpha/\Q_p)^G$,
and if $n-1$ is odd, then $d'$ is even by the Hodge symmetry).
Put
\begin{align*}
e&=\omega\wedge D\omega\wedge\cdots\wedge D^{d'-1}\omega
\in \det[H^{n-1}_\dR(X_{\Q_p}/S_{\Q_p})^G],\\
\wh e&=\wh\omega(1)\wedge\cdots\wedge\wh\omega(s)\wedge
\wh\omega_{s+1}\wedge\cdots\wedge \wh\omega_{d'}\in\Q_p((z))\ot \det[H^{n-1}_\dR(X_{\Q_p}/S_{\Q_p})^G].
\end{align*}
Put $A=\sum_{i=1}^{d'}a_i$ and $B=\sum_{j=1}^{d'}b_j=s+\sum_{j>s}b_j$.
Then $De=(d'-B+Az)/(1-z)e$ as
$P(\ul a';\ul b')=(1-z)D^{d'}+(B-d'-Az)D^{d'-1}+\cdots$,
and
$D\wh e=(\sum_{k>s}(1-b_k))\wh e
=(d'-B)\wh e$ by Lemmas \ref{wh.omega.lem} and \ref{base.lem3}.
Therefore
$z^{B-d'}(1-z)^{d'-A-B}e$ and 
$z^{B-d'}\wh e$ are the horizontal sections,
in particular, they agree up to $\Q_p^\times$.
There is a constant $q\in \Q_p^\times$ such that
\[
\Phi_{\rig,\tau}(z^{B-d'}(1-z)^{d'-A-B}e)=qz^{B-d'}(1-z)^{d'-A-B}e
\iff
\Phi_{\rig,\tau}(e)=q
\overbrace{\left(\frac{z^{B-d'}(1-z)^{d'-A-B}}{(z^\tau)^{B-d'}
(1-z^\tau)^{d'-A-B}}\right)}^{u(z)}e.
\]
Since $u(z)^m$ is a rational function for $m>0$ such that $mA,mB\in\Z$, 
the value $u(\alpha)$ is a $m$-th root of unity.
Therefore we have
$q\sim p^{\frac12(n-1)d'}$, 
where $x\sim y$ means $x=\zeta y$ with some root $\zeta\in \Q_p^\times$ of unity.
On the other hand, with use of the basis $z^{B-d'}\wh e$, one has
\[
\Phi_{\rig,\tau}(z^{B-d'}\wh e)=qz^{B-d'}\wh e\iff \Phi_{\rig,\tau}(\wh e)
=qc^{d'-B}z^{(B-d')(1-p)}\wh e.
\]
Notice that
$B-d'=\sum_{k>s}(b_k-1)=
\sum_{k>s}(b_k^{(1)}-1)$ and $\sum_{k>s}\mu_k=(B-d')(1-p)$.
It follows from \eqref{Katz-thm-eq1-pre} and \eqref{Katz-thm-eq2-pre} that
we have
\[
\Phi_{\rig,\tau}(\wh e)=\left(\ve^{d'}\cdot p^{1+2+\cdots+(s-1)}\,
\prod_{k>s}(\gamma'_k/\gamma'_1)\right) c^{d'-B}z^{(B-d')(1-p)}\wh e,
\]
and hence $q=\ve^{d'}\, p^{\frac12s(s-1)}\prod_{k>s}(\gamma'_k/\gamma'_1)$.
Since $q\sim p^{\frac12(s-1)d'}$ (note $n=s$), it turns out that
\begin{equation}\label{Katz-thm-pf-eq1}
\ve^{d'}\,\prod_{k>s}(\gamma'_k/\gamma'_1)\sim p^{\frac12(s-1)(d'-s)}.
\end{equation}
We compute $\prod_{k>s}(\gamma'_k/\gamma'_1)$.
Since $(b_1,\ldots,b_{d'})$ agrees with
$(b^{(1)}_1,\ldots,b^{(1)}_{d'})$ up to permutation, 
one has
\[
\gamma_1'=p^{-s}\prod_{i=1}^{d'}\frac{\Gamma_p(b_i)}{\Gamma_p(a_i)},\quad
\prod_{k>s}\gamma'_k=\prod_{k>s}\gamma_k
\] by definition.
Notice that the sets $\{a_1,\ldots, a_{d'}\}$ and $\{b_{s+1},\ldots,b_{d'}\}$
are stable under the action $x\mapsto 1-x$. Hence, 
$\prod_{i=1}^{d'}{\Gamma_p(b_i)}/{\Gamma_p(a_i)}$
 is at most a $4$-th root of unity by
 the reflection formula on $p$-adic gamma function
(e.g. \cite[II. 6. (2)]{Ko}, \cite[11.6.12]{Cohen2}).
Therefore we have $\gamma'_1\sim p^{-s}$.
On the other hand,
\begin{align*}
\prod_{k>s}\gamma_k
&=\prod_{k>s}\left((-1)^{Z(b_k)}p^{-Z^\FS(b^{(1)}_k)}\frac{K^\FS(b^{(1)}_k)}{K(b_k)}
\prod_{i=1}^{d'}\frac{\Gamma_p(\{b_i-b_k\})}{\Gamma_p(\{a_i-b_k\})}\right)\\
&=\prod_{k>s}(-1)^{Z(b_k)}p^{-Z^\FS(b^{(1)}_k)}\cdot
\prod_{i=1}^{d'}\prod_{k>s}\frac{\Gamma_p(\{b_i-b_k\})}{\Gamma_p(\{a_i-b_k\})}\\
&=(-p)^{-\delta}
\prod_{i=1}^{d'}\prod_{k>s}\frac{\Gamma_p(\{b_i-b_k\})}{\Gamma_p(\{a_i-b_k\})},
\qquad \left(\delta:=\sum_{k>s}Z(b_k)=\sum_{k>s}Z^\FS(b^{(1)}_k)\right),
\end{align*}
and the last product is again a $4$-th root of unity by
the reflection formula on $p$-adic gamma function.
Therefore we have
\[
\prod_{k>s}\gamma'_k/\gamma'_1=
(\gamma'_1)^{-(d'-s)}\prod_{k>s}\gamma'_k\sim p^{-\delta+s(d'-s)}.
\]
Applying this to \eqref{Katz-thm-pf-eq1}, we have
\[
\ve^{d'}\sim
p^{-\frac12 (d'-s)(s+1)+\delta}.
\] 
We show that the right hand side is equal to $1$, namely
\begin{equation}\label{Katz-thm-eq4}
\sum_{k>s}Z(b_k)=\frac12 (d'-s)(s+1).
\end{equation}
By definition 
\[
Z(x)=\sharp\{i\in\{1,2,\ldots,d'\}\mid a_i<x\}-\sharp\{j\in\{s+1,s+2,\ldots,d'\}\mid b_j<x\}
\]
for $0<x<1$.
Since
$\sum_{k>s}\sharp\{j\in\{s+1,\ldots,d'\}\mid b_j<b_k\}=1+2+\cdots+(d'-s-1)=(d'-s)(d'-s-1)/2$,
\eqref{Katz-thm-eq4} is equivalent to
\[
\sum_{k>s}\sharp\{i\in\{1,2,\ldots,d'\}\mid a_i<b_k\}=\frac12 (d'-s)d'.
\]
Since the set $\{b_{s+1},\ldots,b_{d'}\}$ is stable under the action $x\mapsto 1-x$,
the left hand side can be replaced with
\[
\frac12\sum_{k>s}\bigg(\sharp\{i\in\{1,2,\ldots,d'\}\mid a_i<b_k\}
+\sharp\{i\in\{1,2,\ldots,d'\}\mid a_i<1-b_k\}\bigg).
\]
Now it is immediate to see that the inside of the bracket is equal to $d'$ for each $k$, 
so the summation is $d'(d'-s)/2$ as required.
This completes the proof of \eqref{Katz-thm-eq4}, and hence $\ve$ is a root of unity.

\medskip

There remains to show that 
$\Psi_m^{[s]}(\ul a',\ul b')=\Psi^{[s]}_m(\ul a,\ul b)$.
Recall from Definition \ref{Psi-defn} that
\[
\psi^{(r)}_p(\ul a',\ul b')=\sum_i
\wt\psi^{(r)}_p(a'_i)-\wt\psi^{(r)}_p(b'_i),
\]
\[
\left(\prod_{j>s}\frac{1-x \beta^\FS_j}{1-x \beta_j}\right)
\exp\left(\sum_{r=1}^\infty \psi^{(r-1)}_p(\ul a',\ul b')\frac{x^r}r\right)
=\sum_{m=0}^\infty \Psi_m^{[s]}(\ul a',\ul b')x^m
\]
where $ \beta_j:=1/(b_j-1)$ and 
$\beta^\FS_j:=p^{-1}/(b_j^{(1)}-1)$, and similarly for $\Psi^{[s]}_m(\ul a,\ul b)$.
Therefore it is immediate to see $\Psi_m^{[s]}(\ul a',\ul b')=\Psi^{[s]}_m(\ul a,\ul b)$.
The last assertion follows from Theorem \ref{zeta-value}, namely
\begin{align*}
\psi^{(r)}_p(\ul a,\ul b)&=(d-d^{r+1})L_p(r+1,\omega^{-r})
-\sum_{w_i>1}(w_i-w_i^{r+1})L_p(r+1,\omega^{-r})\\
&=\left((d-d^{r+1})
-\sum_{w_i>1}(w_i-w_i^{r+1})\right)\frac{p^{r+1}-1}{p^{r+1}}\zeta_p(r+1),\quad (r\ne 0),\\
\psi^{(0)}_p(\ul a,\ul b)&=-dp^{-1}\log(d^{p-1})+\sum_{w_i>1}w_ip^{-1}\log(w_i^{p-1}).
\end{align*}
\end{proof}
\begin{rem}\label{Katz-remark}
The proof of Theorem \ref{Katz-thm} shows that $\ve^l=1$ where $l:=\mathrm{lcm}
(4,d',d,w_0,\ldots,w_n)$.
\end{rem}
\subsection{Frobenius on log crystalline cohomology}
\label{log-crys-sect}
We refer to \cite{Ka-log} for a general reference of log schemes and log crystalline cohomology.
For a regular scheme $V$ and a divisor $D$ in $V$, let $(V,D)$ denote the log scheme endowed with
the log structure defined by the reduced part $D_{\mathrm{red}}$ of the divisor $D$.

\bigskip
Let $W=W(k)$ be the Witt ring of a perfect field $k$ of characteristic $p>0$,
and $K=\Frac( W)$ the fractional field. 
Let $\cX$ be a regular $W$-scheme, and let
 \[\xymatrix{f:\cX\ar[r]& \Delta:=\Spec W[[z]]}\]
 be a projective flat $W$-morphism that satisfies the following conditions.
 \begin{itemize}
\item[(i)] Put $O:=\Spec W[[z]]/(z)$. Then $f$ is smooth over $\Delta\setminus O$.
\item[(ii)] Put $Y:=f^{-1}(O)$.
Then the reduced part $Y_{\mathrm{red}}$ 
is a relative simple normal crossing divisor, and
the multiplicity
of each component is prime to $p$.
\end{itemize}
Then, the morphism
\[
\xymatrix{
f:(\cX,Y)\ar[r]&(\Delta,O)
}\]
is log smooth (\cite[Proposition (3.12)]{Ka-log}).
Let $\sigma$ be the $p$-th Frobenius on $W[[z]]$ given by $\tau(z)=cz^p$ with $c\in 1+pW$.
Write $\cX_k=\cX\times_W\Spec k$ etc.
Let
\[
H^j_{\text{\rm log-crys}}(\cX_k/\Delta):=
H^j_{\text{log-crys}}((\cX_{k},Y_{k})/(\Delta,O))
\]
be the $j$-th log crystalline cohomology group
of $(\cX_{k},Y_{k})/(\Delta,O)$ (\cite[\S 6]{Ka-log}).
It is endowed with the $\tau$-linear $p$-th Frobenius, which we write by $\Phi_\crys$.
Letting
$\omega^\bullet_\cX:=\Omega^\bullet_{\cX/W}(\log Y)$
be the de Rham complex with log pole, 
the de Rham complex of $(\cX,Y)/(\Delta,O)$ is given as follows,
\[
\omega^\bullet_{\cX/\Delta}:=\Coker\left[\frac{dz}{z}\ot\omega^{\bullet-1}_\cX\to \omega^\bullet_\cX\right].
\]
Then it follows from \cite[Theorem (6.4)]{Ka-log} that
there is a canonical isomorphism
\begin{equation}\label{comparison}
H^i_{\text{log-crys}}(\cX_k/\Delta)
\cong H^i_\zar(\cX,\omega_{\cX/\Delta}^\bullet),
\end{equation}
which we refer to as the comparison isomorphism.
Let $i:O\to \Delta$ and $i_Y:Y\to \cX$ be the closed immersions, and let
\[
(Y,L_Y):=i_Y^*(\cX,Y),\quad(Y_k,L_{Y_k}):=i_Y^*(\cX_k,Y_k),\quad (O,L_O):=i^*(\Delta,O)
\]
be the pull-back of the log schemes.
Let 
\[
H^i_{\text{log-crys}}(Y_k/W):=
H^i_{\text{log-crys}}((Y_{k},L_{Y_{k}})/(O,L_O))
\] be the $i$-th log cryctalline
cohomology group. The Frobenius $\tau$ induces an endomorphism on $(O,L_O)$, and
then it induces the $p$-th Frobenius $H^i_{\text{log-crys}}(Y_k/W)$, which we write by the same symbol
$\Phi_\crys$.
This is compatible with $\Phi_\crys$ on $H^i_{\text{log-crys}}(\cX_k/\Delta)$
under the canonical map $H^i_{\text{log-crys}}(\cX_k/\Delta)\to H^i_{\text{log-crys}}(Y_k/W)$.
Let $\omega^\bullet_{Y/O}$ be the de Rham complex of $(Y,L_Y)/(O,L_O)$.
Again, thanks to \cite[Theorem (6.4)]{Ka-log}, 
we have the comparison isomorphism
\begin{equation}\label{comparison2}
H^i_{\text{log-crys}}(Y_k/W)
\cong H^i_{\text{log-dR}}(Y/O):=H^i_\zar(Y,\omega_{Y/O}^\bullet)
\end{equation}
called the Hyodo-Kato isomorphism.
The complex $\omega^\bullet_{Y/O}$ fits into an exact sequence
\[
\xymatrix{
0\ar[r]&\omega_{\cX/\Delta}^\bullet\ar[r]^z
&\omega_{\cX/\Delta}^\bullet\ar[r]&\omega^\bullet_{Y/O}\ar[r]&0.
}
\]
This gives rise to the canonical map
\begin{equation}\label{log-crys-eq1}
\xymatrix{
H^i_\zar(\cX,\omega_{\cX/\Delta}^\bullet)\ot W[[z]]/(z)
\ar[r]&
H^i_\zar(Y,\omega_{Y/O}^\bullet)
}\end{equation}
which can be identified with the canonical map 
$H^i_{\text{log-crys}}(\cX_k/\Delta)\to H^i_{\text{log-crys}}(Y_k/W)$.
\begin{lem}\label{log-crys-lem}
The map \eqref{log-crys-eq1} is injective and the cokernel is a finitely generated torsion $W$-module.
Hence we have a natural isomorphism
$K[[z]]/(z)\ot_{W[[z]]} H^i_{\text{\rm log-crys}}(\cX_k/\Delta)\cong
H^i_{\text{\rm log-crys}}(Y_k/W)\ot\Q$.
\end{lem}
\begin{proof}
The injectivity is immediate from the construction.
Let $N$ be the cokernel of \eqref{log-crys-eq1}, so that there is an exact sequence
\[
\xymatrix{
0\ar[r]&N\ar[r]&H^{i+1}_\zar(\cX,\omega_{\cX/\Delta}^\bullet)\ar[r]^z&
H^{i+1}_\zar(\cX,\omega_{\cX/\Delta}^\bullet).}
\]
Since the middle term is a finitely generated $W[[z]]$-module, $N$ is a finitely generated $W$-module.
We want to show $K\ot_WN=0$.
Note that $W[[z]]\to K[[z]]$ is flat as $K[[z]]$ is the $z$-adic completion
of $W[[z]][p^{-1}]$.
Put $\cX_K:=\cX\times_{W[[z]]}K[[z]]$ and $\Delta_K:=\Spec K[[z]]$, and let
$\omega_{\cX_K/\Delta_K}^\bullet$ be the de Rham complex
of $(\cX_K,Y_K)/(\Delta_K,O_K)$.
Then we have an isomorphism
\[
K[[z]]\ot_{W[[z]]}H^{i+1}_\zar(\cX,\omega_{\cX/\Delta}^\bullet)\cong
H^{i+1}_\zar(\cX_K,\omega_{\cX_K/\Delta_K}^\bullet)
\]
as $W[[z]]\to K[[z]]$ is flat, and hence an exact sequence
\[
\xymatrix{
0\ar[r]&K\ot_WN\ar[r]&H^{i+1}_\zar(\cX_K,\omega_{\cX_K/\Delta_K}^\bullet)\ar[r]^z&
H^{i+1}_\zar(\cX_K,\omega_{\cX_K/\Delta_K}^\bullet).}
\] 
Now the vanishing of $K\ot_WN$ follows from the fact that
$H^{j}_\zar(\cX_K,\omega_{\cX_K/\Delta_K}^\bullet)$ is a free $K[[z]]$-module of finite rank
(\cite[(2.18)]{steenbrink}).
\end{proof}
For a sub $K((z))$-module $M\subset 
K((z))\ot H^i_\zar(\cX_K,\omega_{\cX_K/\Delta_K}^\bullet)$
stable under the action of the differential operator $D=z\frac{d}{dz}$, 
we put
$M_{K[[z]]}:=M\cap H^i_\zar(\cX_K,\omega_{\cX_K/\Delta_K}^\bullet)$.
This is Deligne's canonical extension, namely it is
the unique sub $K[[z]]$-module of $M$ of finite rank such 
that $K((z))\ot M_{K[[z]]}=M$, $DM_{K[[z]]}\subset M_{K[[z]]}$ and
that every eigenvalues of the
$K$-linear endomorphism on $M_{K[[z]]}/zM_{K[[z]]}$
given by $v\mapsto Dv$
lie in $[0,1)$.
This can be checked from the fact that 
$H^i_\zar(\cX_K,\omega_{\cX_K/\Delta_K}^\bullet)$ is the canonical extension 
(\cite[(2.20)]{steenbrink}).
We define
\begin{equation}\label{log-crys-eq2}
 \wh M_{\text{\rm log-crys}}
 :=\Image[ M_{K[[z]]}
 \to K[[z]]\ot_{W[[z]]}
H^i_{\text{log-crys}}(\cX_k/\Delta)]\cong M_{K[[z]]}
\end{equation}
a $K[[z]]$-module of finite rank
where the arrow is the comparison map \eqref{comparison}, and
\begin{equation}\label{log-crys-eq3}
 M_{\text{\rm log-crys}}:=\Image[
 \wh M_{\text{\rm log-crys}}
\to  K\ot_W
H^i_{\text{log-crys}}(Y_k/W)]
\end{equation}
 a finite dimensional $K$-module
where the arrow is the canonical map 
$H^i_{\text{log-crys}}(\cX_k/\Delta)\to H^i_{\text{log-crys}}(Y_k/W)$.
Thanks to the property of Deligne's canonical extension (e.g. \cite[III \S 4]{Meb}), 
the arrow $M_{K[[z]]}/zM_{K[[z]]}\to K[[z]]/(z)\ot
H^i_\zar(\cX_K,\omega_{\cX_K/\Delta_K}^\bullet)$ is injective.
Hence by Lemma \ref{log-crys-lem} together with \eqref{comparison},
one has a bijection
\begin{equation}\label{whM.modz.M}
\xymatrix{
 K[[z]]/(z)\ot_{W[[z]]}  \wh M_{\text{\rm log-crys}}\ar[r]^-\cong
 & M_{\text{\rm log-crys}}.
 }\end{equation}

\paragraph{Frobenius on log crystalline cohomology of a totally degenerate fiber}
Let $F$ be a number field, and $\Z_F\subset F$ the ring of integers.
Let
\[
\xymatrix{
f_F:X_F\ar[r]& S_F:=\Spec F[z,(z-z^2)^{-1}]
}\]
be a projective smooth morphism of smooth $F$-schemes.
Let $\mathcal Q$ be a finite set of primes of $F$ 
such that 
$R:=\mathcal Q^{-1}\Z_F$ is etale over $\Z$ and that there is an integral model $f_R:X_R\to S_R$ over $R$
which fits into a commutative diagram
\[
\xymatrix{
Y:=f^{-1}(O)\ar[r]\ar[d]&\cX\ar[d]_{f}&\ar[l]\cX^*\ar[r]\ar[d]&X_{R}\ar[d]^{f_R}\\
O:=\Spec R[[z]]/(z)\ar[r]&\Spec R[[z]]&\ar[l]\Spec R[[z]][z^{-1}]\ar[r]&S_{R}
}
\]
where every squares are cartesian.
We assume that 
$Y$ is reduced and 
a relative simple normal crossing divisor,
which means 
that every components $Y_i$
are smooth and geometrically connected over $R$, and so are every components
of any intersections $Y_{i_1}\cap\cdots\cap Y_{i_r}$.
Let $v$ be a prime in $R$ above $p$, and $k_v$ the residue field.
Let $F_v$ be the $v$-adic completion, and $R_v\subset F_v$ the valuation ring.
Let $\tau_v$ be the $p$-th Frobenius on $R_v[[z]]$ given by $\tau_v(z)=cz^p$
with $c\in 1+pR_v$.
Put $\Delta_v:=\Spec R_v[[z]]$, $O_v:=\Spec R_v[[z]]/(z)$ and let $i_v:O_v\to \Delta_v$
be the embedding.
Write $\cX_{R_v}=\cX\times_{R[[z]]}\Spec R_v[[z]]$ and
$\cX_{k_v}=\cX\times_{R[[z]]}\Spec k_v[[z]]$ etc.
Let $(Y_{k_v},L_{Y_{k_v}})/(O_v,L_{O_v})$ be the pull-back of 
$(\cX_{k_v},Y_{k_v})/(\Delta_v,O_v)$ by $i_v$.
Under the above setting, 
we have the log crystalline cohomology groups
\[
H^i_{\text{\rm log-crys}}(\cX_{k_v}/\Delta_v):=
H^i_{\text{\rm log-crys}}((\cX_{k_v},Y_{k_v})/(\Delta_v,O_v)),\quad
H^i_{\text{\rm log-crys}}(Y_{k_v}):=H^i_{\text{\rm log-crys}}((Y_{k_v},L_{Y_{k_v}})/(O_v,L_{O_v})).
\]
They are endowed with
 the $\tau_v$-linear $p$-th Frobenius, which we write by $\Phi_{\crys,\tau_v}$.
The former Frobenius action extends in a natural way on
\[
F_v((z))\ot_{\O(S_F)} H^i_\dR(X_F/S_F)
\]
via the comparison isomorphism \eqref{comparison}, which we write by the same symbol.
The $\tau_v$-linear Frobenius $\Phi_{\rig,\tau_v}$ on the rigid cohomology
\[
H^*_\rig(X_{k_v}/S_{k_v})\cong 
\O(S_{F_v})^\dag\ot_{\O(S_F)} H^i_\dR(X_F/S_F),
\]
 is compatible with $\Phi_{\crys,\tau_v}$
under the embedding into
a bigger ring
\[
R_v[[z]]\hookrightarrow \left(\varprojlim_nR_v/p^n[[z]][z^{-1}]\right)\ot\Q \hookleftarrow\O(S_{F_v})^\dag.
\]
Let $\cD_F=F[z,(z-z^2)^{-1},\frac{d}{dz}]$, 
and write $H_{P,F}=
\cD_F/\cD_FP$ for $P\in \cD_F$.
For $M_{S_F}\subset H_\dR^i(X_F/S_F)$ a sub $\cD_F$-module,
let $M_{S_F,F[[z]]}$ denote Deligne's canonical extension (cf. \S \ref{base.sect}).
Moreover we put $M_v:=F_v((z))\ot_{\O(S_F)} M_{S_F}$, and 
\[
\wh M_{\text{\rm log-crys},v}:=\wh {(M_v)}_{\text{\rm log-crys}}
\subset
F_v[[z]]
\ot_{R_v[[z]]} H^i_{\text{\rm log-crys}}(\cX_{k_v}/\Delta_v),
\]
\[
M_{\text{\rm log-crys},v}:={(M_v)}_{\text{\rm log-crys}}
\subset
H^i_{\text{\rm log-crys}}(Y_{k_v}/R_v)\ot\Q
\]
the submodules \eqref{log-crys-eq2}, \eqref{log-crys-eq3} arising from $M_v$.

\begin{thm}\label{degenerate.thm}
Let $n\geq 1$.
Let $M_{S_F}\subset H_\dR^{n-1}(X_F/S_F)$ be a sub $\cD_F$-module such that 
$M_v$ is stable under the action of the Frobenius $\Phi_{\crys,\tau_v}$.
Suppose that
there is a finite subset $I\subset(\Z_{(p)}\cap (0,1))^n$ such that there is
an isomorphism
 \begin{equation}\label{degenerate.thm.eq1}
M_{S_F}\cong \bigoplus_{\ul a\in I}  H_{P(\ul a;\ul 1),F}
 \end{equation}
 of left $\cD_F$-modules, where $\ul 1=(1,\ldots,1)$.
 Then, there is
 a $F[[z]]$-basis $\{\wh\omega_{\ul a}(i)\mid\ul a\in I,1\leq i\leq n\}$ of $M_{S_F,F[[z]]}$ such that
 \begin{equation}\label{degenerate.thm.eq3}
p^{-m+1}\Phi_{\crys,\tau_v}(\wh\omega_{\ul a^{(1)}}(m))=\sum_{k=1}^m\left(
\sum_{i+j=m-k}\frac{(-p^{-1}\log c)^i}{i!}
\Psi_j(\ul a,\ul 1)\right)\wh\omega_{\ul a}(k)\,\text{ in }\wh M_{\text{\rm log-crys},v}
\end{equation}
for $1\leq m\leq n$, where we think $\wh\omega(i)$ of being elements of
$\wh M_{\text{\rm log-crys},v}$ by \eqref{log-crys-eq2}.
Moreover, for each $\ul a$, there is an element $\omega_{\ul a}\in M_{S_F}$ which agrees with
$1_\cD+\cD P(\ul a;\ul1)\in H_{P(\ul a;\ul1),F}$ under \eqref{degenerate.thm.eq1} up to $F^\times$,
such that
 \begin{equation}\label{degenerate.thm.eq3.1}
p^{-m+1}\Phi_{\crys,\tau_v}(D^{n-m}\omega_{\ul a^{(1)}})=
\sum_{k=1}^m\left(
\sum_{i+j=m-k}\frac{(-p^{-1}\log c)^i}{i!}
\Psi_j(\ul a,\ul 1)\right)(-1)^{k-m}D^{n-k}\omega_{\ul a}\,\text{ in }
M_{\text{\rm log-crys},v}
\end{equation}
where we think $D^k\omega_{\ul a}$ of being elements of 
$M_{\text{\rm log-crys},v}$ by \eqref{log-crys-eq2} and \eqref{log-crys-eq3}.
\end{thm}
\begin{proof}
Let $\{\wh\omega(i)\}_{i=1,2,\ldots,n}$ be
the $F[[z]]$-basis of 
$\wh H_{P(\ul a;\ul1),F[[z]]}$
in Proposition \ref{base.thm}.
Note that $\wh H_{P(\ul a;\ul1),F[[z]]}$ is Deligne's canonical extension
by Lemma \ref{base.lem4}, so that we have
$M_{S_F,F[[z]]}\cong \bigoplus \wh H_{P(\ul a;\ul1),F[[z]]}$.
Let $\wh\omega_{\ul a}(i)\in M_{S_F,F[[z]]}$ be the associated 
element to $\wh\omega(i)\in \wh H_{P(\ul a;\ul1),F[[z]]}$.
Since $\Phi_{\crys,\tau_v}$ induces the Frobenius intertwiner,
one can apply Theorem \ref{main-c}, and then we have
 \begin{equation}\label{degenerate.thm.eq3-ve}
p^{-m+1}\Phi_{\crys,\tau_v}(\wh\omega_{\ul a^{(1)}}(m))=\ve_{\ul a,v}
\sum_{k=1}^m\left(
\sum_{i+j=m-k}\frac{(-p^{-1}\log c)^i}{i!}
\Psi_j(\ul a,\ul 1)\right)\wh\omega_{\ul a}(k)\,\text{ in }\wh M_{\text{\rm log-crys},v}
\end{equation}
with some constants
 $\ve_{\ul a,v}\in F_v^\times$ depending only on $\ul a$ and $v$.
We shall show that
there is a collection $\{\alpha_{\ul a}\}_{\ul a}$ with $\alpha_{\ul a}\in F^\times$
such that $\ve_{\ul a,v}
=\mathrm{Frob}_v(\alpha_{\ul a^{(1)}})/\alpha_{\ul a}$ for any $v$, where
$\mathrm{Frob}_v$ is
the $p$-th Frobenius on $F_v$ (note that $F_v/\Q_p$ is unramified).
Then, replace $\wh\omega_{\ul a}(i)$ with $\alpha_{\ul a}^{-1}\wh\omega_{\ul a}(k)$,
and we obtain the basis that satisfy \eqref{degenerate.thm.eq3}.
Moreover,
taking the reduction of \eqref{degenerate.thm.eq3} modulo 
$z\wh M_{\text{\rm log-crys},v}$, we have
\eqref{degenerate.thm.eq3.1} as
$\wh\omega_{\ul a}(i)\equiv (-1)^{i+1}D^{n-i}
\omega_{\ul a}$ by Lemma \ref{base.lem4}.

We show the existence of the collection $\{\alpha_{\ul a}\}_{\ul a}$.
To do this,
we employ the weight filtration on log crystalline cohomology by Mokrane \cite{Mok}.
The filtration gives rise to natural maps
\[\xymatrix{
\rho:H^{n-1}_{\text{\rm log-crys}}(Y_{k_v}/R_v)\ot\Q
\ar[r]&
H^0_{\text{\rm crys}}(Y^{[n]}_{k_v}/R_v)\ot\Q
}
\]
\[\xymatrix{
\rho:H^{n-1}_{\text{\rm log-dR}}(Y_F/F)
\ar[r]&
H^0_\dR(Y_F^{[n]}/F)
}
\]compatible under the comparison, 
where $Y_{k_v}^{[s]}$ (resp. $Y^{[s]}_F$) is the disjoint union of $s$-fold intersections of the 
components of $Y_{k_v}$ (resp. $Y_F$).
The $\rho$ on log crystalline cohomology satisfies $\rho\circ\Phi_{\crys,v}=p^{n-1}\Phi_{\crys,v}\circ\rho$.
Let $m=n$ in \eqref{degenerate.thm.eq3.1}.
Since $\rho$ annihilates the image of $D$ by the property of monodromy filtration, we have
\begin{equation}\label{Mok-eq1}
\ve_{\ul a,v}\,\rho(\omega_{\ul a})
=\rho\bigg(p^{-n+1}\Phi_{\crys,\tau_v}(\omega_{\ul a^{(1)}})\bigg)
=
\Phi_{\crys,v}\rho(\omega_{\ul a^{(1)}}).
\end{equation}
Let $Z_i$ be the connected components of $Y^{[n]}_F$, and
put $1_{Z_i}=(0,\ldots,1,\ldots,0)\in H^0_\dR(Y^{[n]}_F/F)=\bigoplus_i
\O(Z_i)$.
Let $\alpha_{\ul a,i}\in F$ be defined by
$\rho(\omega_{\ul a})=\sum_i\alpha_{\ul a,i}1_{Z_i}$.
Then the equality \eqref{Mok-eq1}
implies that $\mathrm{Frob}_v(\alpha_{\ul a^{(1)},i})=\ve_{\ul a,v}\alpha_{\ul a,i}$ in $F_v$
for every $\ul a$ and $i$.
We construct the collection $\{\alpha_{\ul a}\}$ in the following way.
Since the monodromy on $H_{P(\ul a;\ul 1)}$ at $z=0$ is maximal unipotent,
we have $\rho(\omega_{\ul a})\ne0$.
Take an arbitrary $\ul a$ and fix $i$ such that $\alpha_{\ul a,i}\ne0$.
Then $\alpha_{\ul a^{(k)},i}\ne0$ for all $k\geq 0$.
Set
$\alpha_{\ul a^{(k)}}=\alpha_{\ul a^{(k)},i}$.
Let $\ul b\in I\setminus\{\ul a^{(k)}\}_{k\geq0}$ be arbitrary.
In the same way, we set $\alpha_{\ul b^{(k)}}=\alpha_{\ul b^{(k)},j}$ for some $j$.
Continuing this, we obtain $\{\alpha_{\ul a}\}$.
\end{proof}
\begin{thm}\label{degenerate.cor} 
Under the setting in Theorem \ref{degenerate.thm}, let $N>1$ be the minimal integer such that
all $a_i\in 1/N\Z$ for all $\ul a=(a_1,\ldots,a_n)\in I$.
Let $\wt F$ be the Galois closure of $F$ over $\Q$. 
Then, under the Hyodo-Kato isomorphism
\[H^{n-1}_{\text{\rm log-crys}}(Y_{k_v}/R_v)\ot\Q
\cong
 F_v\ot_FH^{n-1}_{\text{\rm log-dR}}(Y_F/F),
\]
the matrix representation of $\Phi_{\crys,v}$ on $M_{\text{\rm log-crys},v}$
with respect to a basis of $H^{n-1}_{\text{\rm log-dR}}(Y_F/F)$ is described by
elements of the $\wt F$-algebra generated by 
\[
\log(c),\quad
\log (m^{p-1}),\quad \chi(m),\quad L_p(r,\chi\omega^{1-r}),
\]
where $m$ and $r$ run over positive integers such that
$m|N$ and $2\leq r\leq n-1$, 
and $\chi$ runs over primitive characters of conductor $f_\chi$ such that
$f_\chi|N$.
\end{thm}
\begin{proof}
This is immediate from
Theorem \ref{degenerate.thm} and Theorem \ref{zeta-value3}.
\end{proof}

\paragraph{Example from \cite{As-Ross}}
Let $d_i>1$ be integers.
Let $U_\Q$ be the affine variety over $S_\Q=\Spec \Q[z,(z-z^2)^{-1}]$
defined by an affine equation
\[
(1-x_1^{d_1})\cdots(1-x_n^{d_n})=z.
\]
We call $U_\Q$ a {\it hypergeometric scheme}.
There is an embedding
\[
\xymatrix{
U_\Q\ar[rr]^-\subset\ar[rd]&&X_\Q\ar[ld]\\
&S_\Q&
}
\]
into a projective smooth scheme $X_\Q$ over $S_\Q$ such that $D=X_\Q\setminus U_\Q$ is a
relative simple normal crossing
divisor over $S_\Q$.
Put $M_{S_\Q}:=\Image[H^{n-1}_\dR(X_\Q/S_\Q)\to H^{n-1}_\dR(U_\Q/S_\Q)]
\cong W_{n-1}H^{n-1}_\dR(U_{\Q}/S_{\Q})$
where $W_\bullet$ denotes the weight filtration.
Let $\mu_m\subset \overline\Q^\times$ denote the group of $m$-th roots of unity.
Let $G=\mu_{d_1}\times\cdots\mu_{d_n}$ act on
$U_{\overline\Q}$ by
\[
(\ve_1,\ldots,\ve_n)\bullet(x_1,\ldots,x_n)
=(\ve_1x_1,\ldots,\ve_nx_n),\quad(\ve_1,\ldots,\ve_n)\in G.
\]  
For $(i_1,\ldots,i_n)\in \Z^{n}$ such that $0<i_k<d_k$, let
$M_{S_\Q}(i_1,\ldots,i_n)$ denote the subspace of $M_{S_\Q}$ on which $(\ve_1,\ldots,\ve_n)$
acts by multiplication by $\ve_1^{i_1}\cdots\ve_n^{i_n}$.
One has
an eigen decomposition
\[
M_{S_\Q}=\bigoplus_{0<i_k<d_k} M_{S_\Q}(i_1,\ldots,i_n).
\]
Put
\begin{equation}\label{As-Ross-eq1}
\omega_{i_1,\ldots,i_n}:=d_1^{-1}x_1^{i_1-d_1}x_2^{i_2-1}\cdots x_n^{i_n-1}
\frac{dx_2\cdots dx_n}{(1-x_2^{d_2})\cdots(1-x_n^{d_n})}\in M_{S_\Q}(i_0,\ldots,i_n).
\end{equation}
Then, 
there is an isomorphism
\[\xymatrix{
 H_{P(\ul a;\ul 1),\Q}\ar[r]^-\cong&
M_{S_\Q}(i_1,\ldots,i_n),\quad Q\longmapsto Q\omega_{i_1,\ldots,i_n}
}
\]
of left $\cD_\Q$-modules,
where $\ul a:=(1-i_1/d_1,\ldots,1-i_n/d_n)$ 
(cf. \cite[Theorems 3.2, 3.3]{As-Ross}).

Fix an embedding $M_{S_\Q}\hra H^{n-1}(X_\Q/S_\Q)$, and we think $M_{S_\Q}$ being a submodule of
$H^{n-1}(X_\Q/S_\Q)$.
One can now apply Theorem \ref{degenerate.thm}.
In this case, one can further show that $\omega_{\ul a}=\omega_{i_1,\ldots,i_n}$ satisfies
\eqref{degenerate.thm.eq3.1}
by the same argument as
the proof of Theorem \ref{degenerate.thm} 
with employing \cite[Theorem 4.6]{As-Ross}.
In particular, when $c=1$, the matrix representation of 
the Frobenius intertwiner
\[
\xymatrix{
\Phi=\Phi_{\rig,\tau}:H_{P(\ul a^{(1)};\ul 1)}^\dag\ar[r]&H_{P(\ul a;\ul 1)}^\dag
}\]
is as follows, 
\[
\begin{pmatrix}
\Phi(\wh\omega_{\ul a^{(1)}}(1))&\cdots&
\Phi(\wh\omega_{\ul a^{(1)}}(n))
\end{pmatrix}=
\begin{pmatrix}
\wh\omega_{\ul a}(1)&\cdots&
\wh\omega_{\ul a}(n)
\end{pmatrix}
\begin{pmatrix}
1&p\Psi_1&p^2\Psi_2&\cdots&p^{n-1}\Psi_{n-1}\\
0&p&p^2\Psi_1&\cdots&p^{n-1}\Psi_{n-2}\\
\vdots&&p^2&\ddots&\vdots\\
&&&\ddots&p^{n-1}\Psi_1\\
0&&\cdots&0&p^{n-1}\\
\end{pmatrix}
\]
where $\Psi_i:=\Psi_i(\ul a;\ul1)$ and 
\begin{align*}
\Psi_1&=\sum_{i=1}^n\wt\psi^{(0)}_p(a_i)=\text{(a linear combination of $\log(m)$'s)}\\
\Psi_2&=\frac12\Psi_1^2+\frac12
\sum_{i=1}^n\wt\psi^{(1)}_p(a_i)
=\frac12\Psi_1^2+\text{(a linear combination of $L_p(2,\chi\omega^{-1})$'s)}\\
\Psi_3&=-\frac13\Psi_1^3+\Psi_1\Psi_2+\frac13
\sum_{i=1}^n\wt\psi^{(2)}_p(a_i)
=-\frac13\Psi_1^3+\Psi_1\Psi_2+\text{(a linear combination of $L_p(3,\chi\omega^{-2})$'s)}\\
&\vdots
\end{align*}
see Theorem \ref{zeta-value3} for the explicit descriptions of the linear combinations of
$p$-adic $L$-values.

\section{Applications to syntomic regulators}\label{syn.sect}
Let $W=W(k)$ be the Witt ring of a perfect field $k$ of charcteristic $p>0$,
and $K=\Frac( W)$ the fractional field. 
We use the notation in \S \ref{notation.sect}.
Write $H_{P,K}=\cD_K/\cD_K P$ for $P\in \cD_K$.
\begin{lem}\label{exist.lem}
Let $(\ul a,\ul b_0)$ and $(\ul a^{(1)},\ul b_0^{(1)})$ be as in Theorem \ref{main.formula1}.
Let $\tau(z)=cz^p$ with $c\in 1+pW$.
Then there is a $\tau$-linear Frobenius intertwiner
\[
\xymatrix{\Phi_\tau:
\whO\ot_{\O_K} H_{P(\ul a^{(1)};\ul b^{(1)}_0)D,K}
\ar[r]& \whO\ot_{\O_K} H_{P(\ul a;\ul b_0)D,K}
}\] 
such that $\Phi_\tau(\whO\ot H_{P(\ul a^{(1)};\ul b^{(1)}_0),K})
\subset \whO\ot H_{P(\ul a;\ul b_0),K}$, where
the inclusion $H_{P(\ul a;\ul b_0),K}\to H_{P(\ul a;\ul b_0)D,K}$ is given by $\theta\mapsto \theta D$ for $\theta\in\cD_K$.
\end{lem}
\begin{proof}
For $\alpha\in \Z_{(p)}$ with $0<\alpha<1$, we put $
Q_\alpha:=P(a_1,\ldots,a_n,\alpha;\overbrace{1,\ldots,1}^{s+1},b_{s+1},\ldots,b_n)$
the hypergeometric differential operator of rank $n+1$.
Let $Q_\alpha^\FS$ be defined in the same way for $\ul a^{(1)},\ul b_0^{(1)}$ and $\alpha^{(1)}$.
Suppose $K=\Q_p$ and $c=1$.
We have the $\tau$-linear Frobenius
\[
\xymatrix{
\Phi_\alpha:H^\dag_{Q_\alpha^\FS}\ar[r]&H^\dag_{Q_\alpha}
}\]
constructed in \S \ref{Frobenius.sect}.
Taking the limit $\alpha\to0$, we have the required
Frobenius intertwiner $\Phi_\tau$
by Lemma \ref{cont-lem}.
For general $K$ and $c$,
the $\tau$-linear $\Phi_\tau$ is obtained by the transformation formula \eqref{change-Frob}.
\end{proof}

\begin{lem}\label{unique.lem}
The $\tau$-linear Frobenius intertwiner in Lemma \ref{exist.lem}
is unique up to multiplication by a nonzero constant.
\end{lem}
\begin{proof}
In this proof, we omit to write ``$K$'' in $H_{Q,K}$ for simplicity. 
One can not directly apply the main theorem of \cite{Dwork-unique}, since
$H_{P(\ul a;\ul b_0)D}$ is not irreducible.
Let $\Phi_i:\tau^*H_{P(\ul a^{(1)};\ul b^{(1)}_0)D}\os{\cong}{\to}H_{P(\ul a;\ul b_0)D}$ ($i=1,2$)
be two $\tau$-linear Frobenius intertwiners.
The map $\Phi_1^{-1}\circ\Phi_2$ is a $\whO$ and $\cD$-linear endomorphism 
of $\whO\ot H_{P(\ul a;\ul b_0)D}$.
In the proof of \cite[Lemma]{Dwork-unique}, Dwork shows that 
if the conditions $(\alpha)$, $(\beta_1)$,\ldots, $(\beta_3)$ in loc.cit. are satisfied, then
any $\whO$ and $\cD$-linear endomorphism 
of $\whO\ot_\O H_{P(\ul a;\ul b_0)D}$ is an algebraic $\cD$-linear endomorphism of
$H_{P(\ul a;\ul b_0)D}$.
Therefore it is enough to show $\End_\cD(H_{P(\ul a;\ul b_0)D})=K$.
An exact sequence
\begin{equation}\label{unique.eq1}
\xymatrix{
0\ar[r]&H_{P(\ul a;\ul b_0)}\ar[r]&H_{ P(\ul a;\ul b_0)D}\ar[r]&\cD/\cD D\ar[r]&0}
\end{equation}
induces an exact sequence
\[
\xymatrix{
0\ar[r]&\Hom_\cD(H_{P(\ul a;\ul b_0)D},H_{P(\ul a;\ul b_0)})\ar[r]&
\End_\cD(H_{P(\ul a;\ul b_0)D})\ar[r]
&\Hom_\cD(H_{P(\ul a;\ul b_0)D},\cD/\cD D).
}\]
There is no left $\cD$-linear map $H_{P(\ul a;\ul b_0)}\to\cD/\cD D$ other than zero.
Indeed, if there is such a map, it must be an isomorphism as
the both sides are irreducible.
When the rank $n$ of $H_{P(\ul a;\ul b_0)}$ is $>1$, this is obviously impossible.
When $n=1$, then $P(\ul a;\ul b_0)=(1-z)D-a_1z$ with $0<a_1<1$ and therefore
it has a non-trivial monodromy at $z=1$. 
Hence
$H_{P(\ul a;\ul b_0)}$ cannot be isomorphic to $\cD/\cD D$.
We now have $\Hom_\cD(H_{P(\ul a;\ul b_0)},\cD/\cD D)=0$ and hence
\[
\xymatrix{
K=\End_\cD(\cD/\cD D)\ar[r]^-\cong&\Hom_\cD(H_{P(\ul a;\ul b_0)D},\cD/\cD D)
}
\]
by the exact sequence \eqref{unique.eq1}.
The rest is to show 
$\Hom_\cD(H_{P(\ul a;\ul b_0)D},H_{P(\ul a;\ul b_0)})=0$.
Suppose that there is a non-zero left $\cD$-linear map
$\psi:H_{P(\ul a;\ul b_0)D}\to H_{P(\ul a;\ul b_0)}$.
This must to be surjective as $H_{P(\ul a;\ul b_0)}$ is irreducible.
Since $\End_\cD(H_{P(\ul a;\ul b_0)})=K$, the composition
$H_{P(\ul a;\ul b_0)}\to H_{P(\ul a;\ul b_0)D}\os{\psi}{\to} H_{P(\ul a;\ul b_0)}$
is multiplication by a nonzero constant.
This means that the exact sequence \eqref{unique.eq1} splits.
Therefore there is a lifting 
$1_\cD+\theta D\in H_{P(\ul a;\ul b_0)D}$ with $\theta\in H_{P(\ul a;\ul b_0)}$
that satisfies $D(1_\cD+\theta D)=0$
in $H_{P(\ul a;\ul b_0)D}$ or equivalently $1_\cD+D\theta=0$ in $H_{P(\ul a;\ul b_0)}$.
We show that this is impossible, namely
\[
1_\cD\not\in \Image[D:H_{P(\ul a;\ul b_0)}\to H_{P(\ul a;\ul b_0)}].
\]
Suppose $1_\cD\in DH_{P(\ul a;\ul 1)}$. 
Let $\O=K[[z]]$ and put $H_e=\sum_{i\geq0}\O D^i$.
Then there is some $N\geq 0$ such that
$1_\cD\in D(z^{-N}H_e)$. 
Let $\overline D$ be the linear endomorphism of $H_e/zH_e$ induced from 
the left action of $D$.
The eigenvalues of $\overline D$ are $0$ and $1-b_j$ ($j>s$).
Therefore, the left action of $D$ on $z^{-i}H_e/z^{-i+1}H_e$ is bijective for all $i\geq1$,
and hence, by the inductive argument, one has $1_\cD\in D(H_e)$.
This implies that $\overline D$ is surjective (hence bijective).
However, this is impossible as $\overline D$ has an eigenvalue $0$
(note $s>0$ by the assumption in Theorem \ref{main.formula1}).
This completes the proof. 
\end{proof}

Summing up the results in Lemmas \ref{unique.lem} and \ref{exist.lem}
and Theorem \ref{main-c},
we have the following.
\begin{thm}\label{syn-thm}
Let the notation and the setting be as in Theorem \ref{main.formula1}.
Let $\tau(z)=cz^p$ with $c\in 1+pW$.
Then a $\tau$-linear Frobenius intertwiner
\[
\xymatrix{\Phi_\tau:
\whO\ot_{\O_K} H_{P(\ul a^{(1)};\ul b^{(1)}_0)D,K}
\ar[r]& \whO\ot_{\O_K} H_{P(\ul a;\ul b_0)D,K}
}\] 
exists, and it is unique up to multiplication by a nonzero constant.
Let $e\in H_{P(\ul a^{(1)};\ul b^{(1)}_0)D,K}$
be the element represented by $1_\cD$. 
Let us choose $\Phi_\tau$ such that 
$\Phi_\tau(e)\equiv e$ modulo the image of 
$\whO\ot_{\O_K} H_{P(\ul a;\ul b_0),K}$.
Then
\[
\Phi_\tau(e)-e=
\sum_{k=1}^{s}\left(
\sum_{i+j=s+1-k}\frac{(-p^{-1}\log c)^i}{i!}
\Psi_j^{[s+1]}(\ul a,\ul b_0)\right)\wh\omega(k)\in \whO\ot H_{P(\ul a;\ul b_0),K}.
\]
\end{thm}

\medskip

We give two examples of application of Theorem \ref{syn-thm} with
$s=n$, i.e. $\ul b_0=(1,\ldots,1)$.

\begin{exmp}[Higher Ross symbols, \cite{As-Ross}]
Let $W=W(k)$ be the Witt ring of an algebraically close field $k$ of characteristic $p$.
Let $K=\Frac(W)$ and $\cD_K:=K\ot_{\Q_p}\cD$.
Recall from \S \ref{Katz} the hypergeometric scheme $U/S$ over $W$
defined by an equation
\[
(1-x_0^{d_0})\cdots(1-x_n^{d_n})=z.
\]
For $d_i$-th roots $\nu_i$ of unity, we consider a Milnor symbol
\[
\xi:=\left\{\frac{1-x_0}{1-\nu_0 x_0},\ldots,\frac{1-x_n}{1-\nu_n x_n}\right\}
\in K^M_{n+1}(\O(U)),\]
 and call it a {\it higher Ross symbol} (loc.cit.).
 One directly has
 \begin{equation}\label{Higher-Ross-eq1}
\dlog(\xi)=(-1)^n\sum_{i_0=1}^{d_0-1}\cdots
\sum_{i_n=1}^{d_n-1}(1-\nu_0^{i_0})\cdots(1-\nu_n^{i_n})\,
\omega_{i_0\ldots i_n}\frac{dz}{z},
\end{equation}
where $\omega_{i_0\ldots i_n}$ are as in \eqref{As-Ross-eq1}.
In \cite[\S 2]{AM}, we construct a symbol map
\begin{equation}\label{AM-symbol}
\xymatrix{
[-]_{U/S}:K^M_{n+1}(\O(U))\ar[r]&
\Ext^1_{\FilFMIC(S)}(\O_S,H^{n}(U/S)(n+1))
}\end{equation}
to the group of $1$-extensions in 
the category of filtered $F$-isocrystals,
where $H^i(U/S)$ denotes the object defined from the
de Rham cohomology $H^i_\dR(U_K/S_K)$ and
the rigid cohomology $H^i_{\mathrm{rig}}(U_{\overline\F_p}/S_{\overline\F_p})$, 
and $\O_S(j)$
denotes the Tate twist.
Therefore, to the higher Ross symbol $\xi$,
the symbol map \eqref{AM-symbol} associates the 1-extension
\begin{equation}\label{Higher-Ross-eq2}
\xymatrix{
0\ar[r]&H^{n}(U/S)\ot\O_S(n+1)\ar[r]&M_\xi\ar[r]&\O_S\ar[r]&0
}
\end{equation}
of filtered $F$-isocrystals on $S$.
Let $e_\xi\in \mathrm{Fil}^0M_\dR$ be the unique lifting of $1\in \O_S(S)$.
Thanks to \cite[Theorem 3.7]{AM},
the element
\[
e_\xi-\Phi_\tau(e_\xi)\in \O(S_K)^\dag \ot H^n_\dR(U_K/S_K)
\] 
agrees with the {\it syntomic regulator} of $\xi$ up to sign in the syntomic cohomology group
of Kato-Tsuji.
A simple argument shows that
$e_\xi-\Phi_\tau(e_\xi)$
lies in $\O(S_K)^\dag \ot W_{n}H^{n}_\dR(U_K/S_K)$, 
so that \eqref{Higher-Ross-eq2} gives rise to
\[
\xymatrix{
0\ar[r]&\cD_K/\cD_K P(\ul a;\ul1)\ar[r]&M_{\ul a}\ar[r]&\O_S\ar[r]&0
}
\]
for each $\ul a=(1-i_0/d_0,\ldots,1-i_n/d_n)$, and the middle term
is isomorphic to $\cD_K/\cD_K P(\ul a;\ul1)D$ by \eqref{Higher-Ross-eq1}.
Therefore one can apply Theorem \ref{syn-thm}, and then
has a full description of the syntomic regulator of $\xi$.
\end{exmp}
\begin{exmp}[{\cite[6.1]{AC}}]
Let $f:Y\to \P^1_{\Z_p}$ be the elliptic fibration whose general fiber
$f^{-1}(z)$ is defined by a Weierstrass equation
\[y^2=x(x^2-2x+1-z).\]
This is smooth over $S=\P^1\setminus\{0,1,\infty\}$.
Let $X=f^{-1}(S)$, and let
\[\xi:=\left\{
\frac{y+x-1}{y-x+1},
\frac{zx}{(1-x)^3}
\right\}
\]
be a Milnor symbol in $K_2$ of $X$.
The symbol map \eqref{AM-symbol} associates a 1-extension
\[
\xymatrix{
0\ar[r]&H^1(X/S)\ot\O_S(2)\ar[r]&M_\xi\ar[r]&\O_S\ar[r]&0
}
\]
of filtered $F$-isocrystals on $S$.
Let $e_\xi\in \mathrm{Fil}^0M_\dR$ be the unique lifting of $1\in \O(S)$.
In \cite[6.1]{AC}, we compute the syntomic regulator
\[
e_\xi-\Phi_\tau(e_\xi)\in H^1_{\mathrm{rig}}(X_{\F_p}/S_{\F_p})\cong \O(S)^\dag\ot
H^1_\dR(X_{\Q_p}/S_{\Q_p})
\] 
explicitly, and provide numerical verifications
of the $p$-adic Beilinson conjecture by Perrin-Riou.
Now, as is computed in \cite[6.1]{AC}, there is an isomorphism
\[M_\dR\cong \cD/\cD P\left(\frac14,\frac34;1,1\right)D\]
of left $\cD$-modules, so that one can apply Theorem \ref{syn-thm}.
The authors of \cite{AC} expect that the constant term of
``$E_2^{(\sigma)}(\lambda)$'' is 
\[
-\frac12\left(\frac1p\log(64^{p-1})+\frac1p\log c\right)^2
\]
(see \cite[6.1, (ix)]{AC}).
This now follows from Theorem \ref{syn-thm} and the fact that
\[
\wt\psi_p^{(0)}\left(\frac14\right)
=\wt\psi_p^{(0)}\left(\frac34\right)=-\frac{1}{p}\log(8^{p-1}),\quad
\wt\psi_p^{(r)}(z)=(-1)^r\wt\psi_p^{(r)}(1-z),
\]
where the former follows from
\cite[Theorem 2.7]{New} and the latter is
\cite[Theorem 2.6 (2)]{New}.
\end{exmp}

\end{document}